\documentclass[final,5p,times]{elsarticle}

\biboptions{numbers,sort&compress}

\usepackage{amsmath,amssymb,amsfonts,amsthm,mathtools}
\usepackage{graphicx}
\usepackage{booktabs}
\usepackage{array}
\usepackage{xcolor}
\usepackage{hyperref}
\hypersetup{pdfauthor={Ismail Belgacem},
  pdftitle={Luenberger observers for two-gene regulatory networks with sigmoidal regulation}}

\newtheorem{theorem}{Theorem}
\newtheorem{lemma}{Lemma}
\newtheorem{proposition}{Proposition}
\newtheorem{corollary}{Corollary}
\newtheorem{assumption}{Assumption}
\newtheorem{remark}{Remark}

\begin{document}

\begin{frontmatter}

\title{Luenberger observers for two-gene regulatory networks with
       sigmoidal regulation: explicit observability conditioning,
       polytopic convergence certificates, and a logistic--Hill
       comparison}

\author{Ismail Belgacem\corref{cor1}}
\cortext[cor1]{Corresponding author.}
\ead{ismail.belgacem.81@gmail.com}
\address{Independent researcher, Mezaourou Ghazaouet,
Tlemcen 13421, Algeria}

\begin{abstract}
Only a few species of a gene regulatory network can be measured, so
the others must be reconstructed by a state observer, and how well
this works depends on the function chosen to describe regulation.  We examine this for the two-gene negative
feedback loop with one measured protein, for monotone regulation
functions including the logistic function and the Hill function
with its leaky and shifted variants.  Using standard tools, we show
that the answer is set by the smallest and the largest slope of
the regulation function on the concentrations visited.  For every function considered, the unmeasured concentration
can be recovered from the measured signal and its rate of change,
with an error sensitivity set by the smallest slope.  A
Luenberger observer with explicit gains converges exponentially from
any initial guess and tolerates bounded measurement noise whenever a
checkable matrix condition involving only the largest slope holds, as
it does for every model examined with bounded slopes.  The guaranteed rate cannot exceed the degradation rate of
the unmeasured protein unless the smallest slope is positive, and faster
rates then hold only for small initial errors.  All
reported rates are certified from above and below.  At matched
threshold and steepness, no family is uniformly better: the logistic
function is better conditioned when the threshold lies mid-range,
leaky or shifted Hill functions when it lies low.  Hill functions with
exponent below one have an unbounded slope at zero concentration and
lose the global guarantee, which logistic and shifted Hill functions
keep.
\end{abstract}


\end{frontmatter}

\section{Introduction}
\label{sec:intro}

In a living cell only a few of the molecular species that make up a
gene regulatory network can be observed, typically through one or two
fluorescent reporters, while the behaviour of the circuit depends on
all of them.  Monitoring a circuit, detecting that it has left its
intended regime, or closing a feedback loop around it therefore
requires an estimate of the concentrations that are not measured.  A
state observer provides such an estimate from a dynamical model of the
network and the available
measurement~\citep{alon2007introduction,dejong2002modeling}.  This
paper studies how the design and the guarantees of such an observer
depend on the way gene regulation is described in the model.

Regulation is usually described by a sigmoidal function of the
regulator concentration.  The most common choice is the Hill
function~\citep{hill1910possible}, often complemented by a basal
(leaky) term or by a shift of its argument so that expression does not
vanish in the absence of the regulator.  The logistic function is
another sigmoid, proposed as a modelling framework for gene networks
in~\citep{belgacem2026mbs}.  The two are closely related, since a Hill
function is exactly a logistic function of the logarithm of the
concentration; they differ in the coordinate in which the sigmoid is
symmetric and in the behaviour of their tails.  Which of them
describes a given regulatory interaction better is an empirical
question, to be settled by dose--response data, and it is not
addressed here.  Our question is narrower and purely mathematical:
given a model of either kind, what can be said about observability and
about the convergence of an observer, and which property of the
regulation function is responsible for each conclusion?

For observers, what matters about a regulation function is its slope.
An upper bound on the slope is what convergence proofs need; a lower
bound on the slope, on the set of concentrations that the network
actually visits, determines how well the unmeasured state can be
distinguished through the measured one.  Both families have explicit
upper bounds in the usual parameter range: the largest logistic slope
is proportional to the steepness parameter and does not depend on the
threshold, and the largest Hill slope is an explicit function of the
Hill exponent and threshold when the exponent is at least one.  In
both cases the bound holds uniformly in the concentration and grows
with the steepness of the sigmoid, so that neither family is
analytically out of reach of standard observer theory.  The families
do differ in the lower bound.  The slope of an unshifted Hill function
with exponent above one vanishes at zero concentration, the logistic
slope does not; on the other hand the logistic slope decays
exponentially at high concentration, where Hill-type slopes decay only
algebraically.  For Hill exponents below one the situation is reversed:
the slope becomes unbounded at zero concentration, and it is the upper
bound that is lost.

We work these effects out for the two-gene negative feedback loop with
one measured protein, the simplest network in which the unmeasured
species is seen only through a regulation function.  The tools are
standard: the Lie-derivative observability
map~\citep{hermann1977nonlinear,besancon2007nonlinear}, the Luenberger
observer~\citep{luenberger1971introduction}, the mean-value (polytopic)
representation of the observer error used for Lipschitz
systems~\citep{rajamani1998observers,zemouche2013lmi}, quadratic
certificates for monotone nonlinearities~\citep{arcak2001nonlinear} and
polytopic matrix inequalities~\citep{boyd1994linear}.  We do not
propose a new observer methodology.  The contribution is an explicit
and complete treatment of one model class in which every statement is
attributed to the property of the regulation function that it actually
uses, so that the comparison between families rests on computed
quantities and not on the choice of a particular Hill variant.

Observers for biological systems have a long history in bioprocess
engineering~\citep{gauthier1992simple}, and the place of control
theory in molecular systems biology is surveyed
in~\citep{sontag2005molecular}.  Recent designs for gene-expression
and bioprocess models include a high-gain observer that reconstructs
the mRNA and protein concentrations of a one-gene circuit with
end-product activation~\citep{bouhadjra2022highgain}, a positive state
observer of Luenberger type for the multidimensional Goodwin
oscillator~\citep{proskurnikov2019positive}, and a high-gain observer
for bilinear systems in block form applied to bioreactor
models~\citep{rezarahmati2023bilinear}.  These works fix the form of
the nonlinearity and exploit its structure.  The question here is
different: how observability and the guarantees of a given observer
depend on the choice of the regulation function itself, among
families that are used interchangeably in practice.

The results are as follows.  For any strictly monotone regulation
function, the measured concentration and its time derivative determine
the unmeasured one, and the sensitivity of this reconstruction is
given in closed form as a function of the slope.  A Luenberger
observer whose gains are given by pole placement has an error equation
that can be written exactly as a linear system with coefficients
confined to a rectangle of slopes; when a matrix inequality at the
four corners of this rectangle is feasible, it gives exponential
convergence from every initial condition, together with a bound on the
effect of measurement noise.  The rate of any such global guarantee is limited
by the degradation rate of the unmeasured protein, a limit that only a
positive lower slope bound can lift.  These statements hold for the
logistic, Hill, leaky Hill and shifted Hill families alike when the
Hill exponent exceeds one, and the inequality is feasible in every
case examined.  Compared on the set to which trajectories
are eventually confined, at matched threshold and steepness, the
unshifted Hill model is observable but with a negligible margin, a
leak or a shift removes most of this degeneracy, and between the
logistic and the generalised Hill forms the ranking depends on where
the threshold lies within the expression range.  Certified convergence
rates, in contrast, differ little between families.

Several things are outside the scope of the paper.  We do not fit
models to data or perform model selection, and we make no claim that
one family is biologically preferable.  The observer assumes that the
model and its parameters are known.  Networks with transmission delays
are not considered.  Control of the same feedback loop under uncertain
measurements was studied
in~\citep{chambon2020qualitative,belgacem2020control}, and the
observer problem treated here is complementary to that work.

The paper is organised as follows.  Section~\ref{sec:model} defines
the model class, the regulation families and the set on which slopes
are evaluated.  Section~\ref{sec:observability} treats observability,
Section~\ref{sec:design} the observer and its convergence
certificates, and Section~\ref{sec:numerics} illustrates them on a
baseline example.  Section~\ref{sec:compare} compares the families for
Hill exponents above one, and Section~\ref{sec:subunit} treats
exponents not exceeding one.  Section~\ref{sec:conclusion} discusses
limitations and concludes, \ref{app:certificates} lists the numerical
certificates behind every reported rate, and \ref{app:chain} extends the
observability computation to a regulatory cascade.  Throughout,
$\|\cdot\|$ is the Euclidean norm and the induced matrix norm, and for
a positive definite matrix $P$ we write
$\kappa(P)=\lambda_{\max}(P)/\lambda_{\min}(P)$, not to be confused
with the rate constants $\kappa_i$.

\section{Model class and operating set}
\label{sec:model}

\subsection{Model}

We consider the two-gene negative feedback loop, the smallest member
of the negative-feedback architectures studied
in~\citep{goodwin1965oscillatory,elowitz2000synthetic},
\begin{equation}
\label{eq:GRN}
\begin{aligned}
  \dot{x}_1 &= \kappa_1 r_1(x_2) - \gamma_1 x_1, \\
  \dot{x}_2 &= \kappa_2 r_2(x_1) - \gamma_2 x_2,
\end{aligned}
\qquad y=x_1,
\end{equation}
with $\kappa_i,\gamma_i>0$: gene~1 is repressed by gene~2 through
$r_1$, gene~2 is activated by gene~1 through $r_2$, and $x_1$ is
measured (the case $y=x_2$ is symmetric).  Loops with two repressions,
such as the toggle switch~\citep{gardner2000construction}, can be
bistable and are not considered.

\begin{assumption}
\label{ass:reg}
$r_1,r_2:\mathbb{R}\to\mathbb{R}$ are $C^1$ and take values in $[0,1]$
on $[0,\infty)$, $r_1$ is non-increasing and $r_2$ non-decreasing, both
are strictly monotone on $[0,\infty)$, and
$0\le -r_1'\le M_1$, $0\le r_2'\le M_2$ on $\mathbb{R}$.
\end{assumption}

\subsection{Regulation families}

Table~\ref{tab:families} lists the activation forms $a(x)$ considered.
The logistic function is $f(x)=1/(1+e^{-\lambda(x-\theta)})$ and the
Hill function is $h(x)=x^n/(K^n+x^n)$; for $x>0$,
\begin{equation}
  \frac{x^n}{K^n+x^n}=\frac{1}{1+e^{-n(\ln x-\ln K)}} ,
  \label{eq:hill_is_logistic}
\end{equation}
so the Hill function is a logistic function of $\ln x$.

The logistic slope is $f'=\lambda f(1-f)\in(0,\lambda/4]$, with
maximum at $x=\theta$.  The Hill slope depends on the exponent.  For
$n>1$ it vanishes at $x=0$ and is bounded by
\begin{equation}
  L_H(n,K)=\frac{(n+1)^2}{4nK}\Bigl(\frac{n-1}{n+1}\Bigr)^{(n-1)/n},
  \label{eq:LH}
\end{equation}
attained at $x=K[(n-1)/(n+1)]^{1/n}$.  For $n=1$ it is at most $1/K$.
For $0<n<1$, $h'(x)\sim nK^{-n}x^{n-1}\to\infty$ as $x\to0^+$: the
function is H\"older continuous but has no finite Lipschitz constant on
an interval containing zero.  Both bounds, $\lambda/4$ and $L_H(n,K)$,
are uniform in the state and depend on the parameters of the function:
the logistic bound on $\lambda$ but not on the threshold $\theta$, the
Hill bound on $n$ and $K$.  Each grows without limit with its steepness
parameter, and for parameters known only within a range a common bound
is obtained in either case by maximising the explicit expression over
that range.  Every argument that needs only an upper slope bound
therefore applies to both families when $n\ge1$.

We set $r_2=a$.  For the logistic, Hill and shifted-Hill families the
repression form is the complement $r_1=1-a$.  For the leaky family we
use
\begin{equation*}
  r_2=b+(1-b)h,\qquad r_1=b+(1-b)(1-h),
\end{equation*}
so that the activation has basal value $r_2(0)=b$ and the repression
has the same leak floor, $\lim_{x\to\infty}r_1(x)=b$; these two
functions are specified separately and are not complements
($r_1\ne1-r_2$).
Both have slope magnitude $(1-b)h'$.  We write $\theta_j$ (or $K_j$)
for the threshold of the regulation exerted by $x_j$.

Hill-type functions are defined for $x\ge0$ only.  The plant never
uses their values elsewhere, since $x_i\ge0$; only the observer state
may leave the positive orthant, and the continuation of $r_i$ to
negative arguments is therefore a choice made in the design of the
observer.  Any $C^1$ monotone continuation that respects the slope
bounds is admissible.  For $n>1$ we extend Hill-type functions by their
value at the left end of their domain ($x=0$, or $x=-x_0$ in the
shifted case), which is $C^1$ because the slope vanishes there, so
Assumption~\ref{ass:reg} holds for all four families.  The logistic
function needs no extension.  Sections~\ref{sec:model}
to~\ref{sec:compare} assume $n>1$ for Hill-type functions;
Section~\ref{sec:subunit} treats $n\le1$.

\begin{table}[ht]
\centering
\caption{Activation forms and their slopes ($n>1$, $0<b<1$, $x_0>0$;
$h$ is the Hill function and $L_H$ is given by~\eqref{eq:LH}).  The
last column is the slope bound on the extended domain $\mathbb{R}$;
for the shifted form the supremum over $x\ge0$ alone is
$\sup_{z\ge x_0}h'(z)$, which is smaller than $L_H$ when the maximiser
of $h'$ lies below $x_0$.}
\label{tab:families}
\small
\begin{tabular}{@{}llll@{}}
\toprule
Family & $a(x)$ & $a'(0)$ & $\sup a'$ \\
\midrule
Logistic     & $1/(1+e^{-\lambda(x-\theta)})$ & $f'(0)>0$ & $\lambda/4$ \\
Hill         & $h(x)=x^n/(K^n+x^n)$           & $0$                    & $L_H(n,K)$ \\
Leaky Hill   & $b+(1-b)h(x)$                  & $0$                    & $(1-b)L_H$ \\
Shifted Hill & $h(x+x_0)$                     & $h'(x_0)>0$            & $L_H(n,K)$ \\
\bottomrule
\end{tabular}
\end{table}

\subsection{Operating set and equilibrium}

The following two statements are standard for this loop; we include
short proofs because the constants are used later.

\begin{lemma}[Invariant and attracting boxes]
\label{lem:boxes}
Let $\mathcal{B}_0$ be the box
$[0,\kappa_1/\gamma_1]\times[0,\kappa_2/\gamma_2]$ and define the boxes $\mathcal{B}_k=[a_1^k,b_1^k]\times[a_2^k,b_2^k]$ by
\begin{equation*}
\begin{aligned}
  a_1^{k+1}&=\tfrac{\kappa_1}{\gamma_1}r_1(b_2^k), &
  b_1^{k+1}&=\tfrac{\kappa_1}{\gamma_1}r_1(a_2^k),\\
  a_2^{k+1}&=\tfrac{\kappa_2}{\gamma_2}r_2(a_1^k), &
  b_2^{k+1}&=\tfrac{\kappa_2}{\gamma_2}r_2(b_1^k).
\end{aligned}
\end{equation*}
Then $\mathcal{B}_{k+1}\subseteq\mathcal{B}_k$, every $\mathcal{B}_k$
is forward invariant for~\eqref{eq:GRN}, and every solution starting
in $\mathcal{B}_0$ is asymptotically confined to
$\mathcal{B}_\infty=\bigcap_k\mathcal{B}_k$, in the sense that
$\textup{dist}(\mathbf{x}(t),\mathcal{B}_\infty)\to0$ as
$t\to\infty$.
\end{lemma}

\begin{proof}
$\mathcal{B}_1\subseteq\mathcal{B}_0$ because $r_i\in[0,1]$, and
nestedness follows by induction from the monotonicity of $r_i$.  On
$\mathcal{B}_{k+1}\subseteq\mathcal{B}_k$ we have
$\kappa_1r_1(b_2^k)-\gamma_1x_1\le\dot x_1\le\kappa_1r_1(a_2^k)-\gamma_1x_1$
and similarly for $x_2$, so on each face of $\mathcal{B}_{k+1}$ the
vector field points inwards or is tangent, and $\mathcal{B}_{k+1}$ is
forward invariant.  For the confinement, let
$\underline x_i=\liminf_{t\to\infty}x_i(t)$ and
$\overline x_i=\limsup_{t\to\infty}x_i(t)$.  The variation-of-constants
formula for each scalar equation and the monotonicity of $r_i$ give
\begin{equation*}
\begin{aligned}
  \underline x_1&\ge\tfrac{\kappa_1}{\gamma_1}r_1(\overline x_2), &
  \overline x_1&\le\tfrac{\kappa_1}{\gamma_1}r_1(\underline x_2),\\
  \underline x_2&\ge\tfrac{\kappa_2}{\gamma_2}r_2(\underline x_1), &
  \overline x_2&\le\tfrac{\kappa_2}{\gamma_2}r_2(\overline x_1).
\end{aligned}
\end{equation*}
Starting from $\underline x_i\ge a_i^0$, $\overline x_i\le b_i^0$ and
iterating gives $\underline x_i\ge a_i^k$ and $\overline x_i\le b_i^k$
for all $k$; the endpoint sequences are monotone and bounded, so they
converge, and the claim follows.
\end{proof}

We call $\mathcal{B}_\infty=J_1\times J_2$ the \emph{operating set}.
It is an invariant enclosure of the long-term behaviour, not
necessarily the exact set of limit points: it need not shrink to the
equilibrium.  The approach is quantified by the inequalities in the
proof: a solution starting in $\mathcal{B}_0$ has
$x_i(t)\ge a_i^{1}-(a_i^{1}-x_i(0))^{+}e^{-\gamma_it}$, so a positive
basal level moves trajectories away from the face $x_i=0$ at the
degradation rate.  The operating set is the natural set on which to
evaluate slope bounds: for instance, a leaky Hill function has zero
slope at $x=0$, but $J_2$ then starts at or above
$\kappa_2b/\gamma_2>0$, where the slope is positive.  On an interval
$J$ we write
\begin{equation}
  m_1(J)=\min_{J}|r_1'|,\qquad m_2(J)=\min_{J}r_2' .
  \label{eq:m_def}
\end{equation}

\begin{proposition}[Equilibrium; classical]
\label{prop:equil}
Under Assumption~\ref{ass:reg}, system~\eqref{eq:GRN} has a unique
equilibrium $\mathbf{x}^*\in\mathcal{B}_\infty$, and it is globally
asymptotically stable on $[0,\infty)^2$.
\end{proposition}

\begin{proof}
Equilibria solve $x_2=G(x_2)$ with
\[
  G(x_2)=\tfrac{\kappa_2}{\gamma_2}\,r_2\bigl(\tfrac{\kappa_1}{\gamma_1}r_1(x_2)\bigr),
\]
which is non-increasing with values in $[0,\kappa_2/\gamma_2]$, so
$x_2-G(x_2)$ is strictly increasing, non-positive at $x_2=0$ and
non-negative at $x_2=\kappa_2/\gamma_2$; it has exactly one zero, and
$\mathbf{x}^*$ lies in every $\mathcal{B}_k$ by induction.  The
Jacobian has trace $-\gamma_1-\gamma_2<0$ and determinant
$\gamma_1\gamma_2+\kappa_1\kappa_2|r_1'|r_2'>0$, so $\mathbf{x}^*$ is
locally exponentially stable.  The divergence is the constant
$-\gamma_1-\gamma_2$, which by the Bendixson criterion excludes
periodic orbits and homoclinic loops; solutions are bounded, so by the
Poincar\'e--Bendixson theorem every solution converges to
$\mathbf{x}^*$.
\end{proof}

For the logistic family the thresholds that realise a prescribed
equilibrium are explicit: with $\beta_i=\gamma_ix_i^*/\kappa_i\in(0,1)$,
$\theta_2=x_2^*+\lambda^{-1}\ln[\beta_1/(1-\beta_1)]$ and
$\theta_1=x_1^*-\lambda^{-1}\ln[\beta_2/(1-\beta_2)]$.  The inversion
is equally elementary for the other families: for the Hill activation
$K_1=x_1^*[(1-\beta_2)/\beta_2]^{1/n}$ and for the Hill repression
$K_2=x_2^*[\beta_1/(1-\beta_1)]^{1/n}$; for the shifted form the same
with $x_j^*$ replaced by $x_j^*+x_0$; for the leaky form the same with
$\beta_i$ replaced by $(\beta_i-b)/(1-b)$, which requires
$b<\beta_i<1$.

Unless stated otherwise, numerical statements use the following baseline
parameter set: the logistic family with
$\kappa_i=1$~nM\,min$^{-1}$, $\gamma_i=0.1$~min$^{-1}$,
$\theta_i=5$~nM and $\lambda=1$~nM$^{-1}$, so that
$\mathbf{x}^*=(5,5)$~nM, $\mathcal{B}_0=[0,10]^2$ and
$\mathcal{B}_\infty\approx[0.072,9.93]^2$.

\section{Observability: injectivity and conditioning}
\label{sec:observability}

Several distinct properties are commonly called observability, and we
keep them apart:
(a)~the \emph{rank condition}, i.e.\ non-singularity of the
Lie-derivative matrix at a point, which gives local weak
observability~\citep{hermann1977nonlinear};
(b)~\emph{global injectivity} of the observability map
$\Phi(\mathbf{x})=(y,\dot y)$, which means that the output and its
derivative determine the state;
(c)~\emph{quantitative conditioning} of $\Phi$ on a prescribed set,
i.e.\ a bound on how much errors in $(y,\dot y)$ are amplified in the
reconstructed state;
(d)~\emph{uniform observability} in the sense used for systems with
inputs and for high-gain
observers~\citep{gauthier1992simple,besancon2007nonlinear}.
System~\eqref{eq:GRN} has no input, so (d) is not at issue; this
section establishes (a)--(c).

The two components of $(y,\dot y)$ have different units (nM and
nM\,min$^{-1}$), so singular values and condition numbers of the
observability map are meaningful only once a time scale is fixed.  We
therefore choose a reference time $\tau>0$ and use
$\Phi(\mathbf{x})=(y,\tau\dot y)
=\bigl(x_1,\;\tau(\kappa_1r_1(x_2)-\gamma_1x_1)\bigr)$, with
\begin{equation}
  \mathcal{O}(\mathbf{x})=\frac{\partial\Phi}{\partial\mathbf{x}}
  =\begin{bmatrix} 1 & 0 \\ -\tau\gamma_1 & \tau\kappa_1 r_1'(x_2)\end{bmatrix}.
  \label{eq:obs_matrix}
\end{equation}
With $c=1+\tau^2\gamma_1^2$, for $d\ge0$ let
\begin{equation}
  s(d)=\Bigl[\tfrac12\Bigl(c+d^2
        -\sqrt{(c+d^2)^2-4d^2}\Bigr)\Bigr]^{1/2}.
  \label{eq:s_def}
\end{equation}
Injectivity and the rank condition do not depend on $\tau$; the
quantitative statements do, and all numerical values below use
$\tau=1$~min.

\begin{theorem}[Observability map of the two-gene loop]
\label{thm:obs}
Let $J\subseteq\mathbb{R}$ be an interval.
\begin{enumerate}
\item[\textup{(i)}] Let $r_1$ be continuous on $J$.  Then $\Phi$ is
  injective on $\mathbb{R}\times J$ if and only if $r_1$ is injective
  on $J$.
\item[\textup{(ii)}] Let $r_1$ be $C^1$.  The rank condition holds at
  $\mathbf{x}$ if and only if $r_1'(x_2)\ne0$.  The singular values of
  $\mathcal{O}(\mathbf{x})$ satisfy
  $\sigma_{\min}\sigma_{\max}=\tau\kappa_1|r_1'(x_2)|$ and
  $\sigma_{\min}^2+\sigma_{\max}^2=c+\tau^2\kappa_1^2r_1'(x_2)^2$,
  so that $\sigma_{\min}=s(\tau\kappa_1|r_1'(x_2)|)$.  The function $s$
  is strictly increasing and
  $d/\sqrt{c+d^2}\le s(d)\le\min\{d,1\}$.
\item[\textup{(iii)}] Let $r_1$ be $C^1$.  If $|r_1'|\ge m_1>0$ on
  $J$, then $\Phi$ is a
  $C^1$ diffeomorphism of $\mathbb{R}\times J$ onto its image and, for
  all $\mathbf{x},\mathbf{x}'\in\mathbb{R}\times J$,
  \begin{equation}
    \|\mathbf{x}-\mathbf{x}'\|\le
    \frac{\|\Phi(\mathbf{x})-\Phi(\mathbf{x}')\|}{s(\tau\kappa_1m_1)} .
    \label{eq:inverse_lip}
  \end{equation}
\end{enumerate}
\end{theorem}

\begin{proof}
(i)~The first component of $\Phi$ is $x_1$; given $x_1$, the second
determines $r_1(x_2)$, hence $x_2$ if $r_1$ is injective on $J$.
Conversely, if $r_1(x_2)=r_1(x_2')$ with $x_2\ne x_2'$ in $J$, then
$\Phi(x_1,x_2)=\Phi(x_1,x_2')$.
(ii)~$\det\mathcal{O}=\tau\kappa_1r_1'(x_2)$ and
$\textup{tr}(\mathcal{O}^T\mathcal{O})=c+\tau^2\kappa_1^2r_1'^2$
give the two relations, and solving the resulting quadratic for
$\sigma_{\min}^2$ gives~\eqref{eq:s_def}.  With $u=d^2$ and
$S=\sqrt{(c+u)^2-4u}$,
$2\,d(s^2)/du=1-(c+u-2)/S>0$ because
$S^2-(c+u-2)^2=4(c-1)=4\tau^2\gamma_1^2>0$.  The lower bound follows from
$\sigma_{\max}\le\|\mathcal{O}\|_F$, and the upper bound from
$\sigma_{\max}\ge\|\mathcal{O}e_1\|=\sqrt c\ge1$ and from the fact that
$\sigma_{\min}$ does not exceed the norm of the first row, which is $1$.
(iii)~For $\mathbf{x},\mathbf{x}'\in\mathbb{R}\times J$,
\begin{equation*}
  \Phi(\mathbf{x})-\Phi(\mathbf{x}')=
  \begin{bmatrix}1&0\\-\tau\gamma_1&\tau\kappa_1\bar c\end{bmatrix}
  (\mathbf{x}-\mathbf{x}'),
\end{equation*}
with the secant slope
$\bar c=\int_0^1r_1'\bigl(x_2'+\vartheta(x_2-x_2')\bigr)d\vartheta$.
Since $r_1'$ does not vanish on $J$ it has constant sign there, so
$|\bar c|\ge m_1$.  The computation in (ii) applies to this matrix: its
smallest singular value is
$s(\tau\kappa_1|\bar c|)\ge s(\tau\kappa_1m_1)$, which
gives~\eqref{eq:inverse_lip}.  By (i) and (ii) $\Phi$ is an injective
local diffeomorphism, hence a diffeomorphism onto its image.
\end{proof}

Part~(i) holds for all four families on $J=[0,\infty)$: the state of a
Hill-based loop is also uniquely determined by $(y,\dot y)$.  Global
injectivity does not imply the rank condition at every point: a
strictly monotone $C^1$ function may have zero slope at isolated
points, the rank condition holds exactly where $r_1'(x_2)\ne0$, and a
uniform conditioning bound on an interval requires
$\inf|r_1'|>0$ there.  The families differ in (ii)--(iii), that is, in
the value of $m_1$ on the interval of interest.

\begin{corollary}[Slope bounds by family]
\label{cor:families}
Let $J=[\underline x,\overline x]$ with $\underline x\ge0$.
\begin{enumerate}
\item[\textup{(a)}] Logistic: with
  $g(x)=f(x)(1-f(x))$, which is unimodal with maximum $1/4$ at
  $\theta$, $m_1(J)=\lambda\min\{g(\underline x),g(\overline x)\}>0$.
  The rank condition holds on all of $\mathbb{R}^2$, but
  $\inf_{\mathbb{R}}f'=0$: there is no positive bound on unbounded
  sets, and $m_1(J)$ decays exponentially as either end of $J$ moves
  away from $\theta$.
\item[\textup{(b)}] Hill, $n>1$: $h'$ is unimodal on $[0,\infty)$, so
  $m_1(J)=\min\{h'(\underline x),h'(\overline x)\}$, which is positive
  if and only if $\underline x>0$.  On $[0,\infty)^2$ the rank
  condition fails only on the line $x_2=0$.
\item[\textup{(c)}] Leaky Hill: as in \textup{(b)} with the factor
  $(1-b)$; on the operating set $\underline x\ge\kappa_2b/\gamma_2>0$.
\item[\textup{(d)}] Shifted Hill: as in \textup{(b)} with $x$ replaced
  by $x+x_0$; $m_1(J)>0$ also for $\underline x=0$, and the rank
  condition holds on $[0,\infty)^2$.
\end{enumerate}
\end{corollary}

For the baseline with $\tau=1$~min, on $J=[0,10]$ we have $m_1=\lambda\mu$ with
$\mu=g(0)=e^{-\lambda\theta}/(1+e^{-\lambda\theta})^2\approx6.65\times10^{-3}$;
$\sigma_{\min}$ ranges from $0.00662$ at the ends of $J$ to $0.249$ at
$x_2=\theta$, and the condition number of $\mathcal{O}$ from $151.9$
to $4.04$.  The margin is positive but small at the ends of the box,
which is the exponential tail decay noted in (a).

\section{Observer design and convergence certificates}
\label{sec:design}

With a measurement $y_m=x_1+w$, where $w:[0,\infty)\to\mathbb{R}$ is
measurable with $|w(t)|\le W$ for almost every $t$ (so that solutions
are absolutely continuous and differential inequalities hold almost
everywhere), we use the Luenberger
observer~\citep{luenberger1971introduction}
\begin{equation}
\begin{aligned}
  \dot{\hat x}_1&=\kappa_1r_1(\hat x_2)-\gamma_1\hat x_1+l_1(y_m-\hat x_1),\\
  \dot{\hat x}_2&=\kappa_2r_2(\hat x_1)-\gamma_2\hat x_2+l_2(y_m-\hat x_1),
\end{aligned}
\label{eq:observer}
\end{equation}
and write $\mathbf{e}=\mathbf{x}-\hat{\mathbf{x}}$,
$\mathbf{L}=[l_1,l_2]^T$.

\subsection{Gains}

Let $a_{12}=\kappa_1r_1'(x_2^*)\le0$, $a_{21}=\kappa_2r_2'(x_1^*)\ge0$
and $p=-a_{12}a_{21}\ge0$.  The linearised error matrix is
$A_{cl}=\bigl(\begin{smallmatrix}-\gamma_1-l_1&a_{12}\\a_{21}-l_2&-\gamma_2\end{smallmatrix}\bigr)$.

\begin{proposition}[Pole-placement gains]
\label{prop:gains}
Assume $a_{12}\ne0$, i.e.\ $r_1'(x_2^*)\ne0$.  Then $A_{cl}$ has
eigenvalues $-\alpha_1,-\alpha_2$ if and only if
\begin{equation}
\begin{aligned}
  l_1&=\alpha_1+\alpha_2-\gamma_1-\gamma_2,\\
  l_2&=\frac{(\alpha_1-\gamma_2)(\alpha_2-\gamma_2)-p}{a_{12}} .
\end{aligned}
  \label{eq:gains}
\end{equation}
If $\alpha_1,\alpha_2>0$ and $|r_1'(x_2^*)|\ge m_1>0$, then
$|l_2|\le[(\alpha_1+\gamma_2)(\alpha_2+\gamma_2)+p]/(\kappa_1m_1)$.
\end{proposition}

\begin{proof}
Match the characteristic polynomial
$\xi^2+(\gamma_1+\gamma_2+l_1)\xi+(\gamma_1+l_1)\gamma_2-a_{12}(a_{21}-l_2)$
with $(\xi+\alpha_1)(\xi+\alpha_2)$.  The bound follows from
$|(\alpha_1-\gamma_2)(\alpha_2-\gamma_2)|\le(\alpha_1+\gamma_2)(\alpha_2+\gamma_2)$
and $|a_{12}|\ge\kappa_1m_1$.
\end{proof}

This is ordinary pole placement for a $2\times2$ pair; we state it
because the formulas are used below.  The hypothesis $a_{12}\ne0$ is
not implied by strict monotonicity; if $a_{12}=0$ the pair is not
observable at the equilibrium and $-\gamma_2$ is a fixed eigenvalue of
$A_{cl}$ for every gain.  The gain $l_2$ grows like
$1/|a_{12}|$, the familiar loss of conditioning as observability
deteriorates, and this is true for every family: the gain is large
whenever the equilibrium sits where the slope of $r_1$ is small.  For
the baseline, $a_{12}=-0.25$, $a_{21}=0.25$, $p=0.0625$, and placing
both poles at $-1$~min$^{-1}$ gives $l_1=1.8$, $l_2=-2.99$ and
$A_{cl}=\bigl(\begin{smallmatrix}-1.9&-0.25\\3.24&-0.1\end{smallmatrix}\bigr)$.

\subsection{Exact error representation}

\begin{lemma}[Mean-value form of the error dynamics]
\label{lem:error}
Let $r_1,r_2$ be $C^1$.  Along any solution of~\eqref{eq:GRN}
and~\eqref{eq:observer},
\begin{equation}
  \dot{\mathbf{e}}=\tilde A(c_1(t),c_2(t))\,\mathbf{e}-\mathbf{L}w,
  \;\;
  \tilde A(c_1,c_2)=
  \begin{bmatrix}
    -(\gamma_1+l_1) & -\kappa_1c_1 \\
    \kappa_2c_2-l_2 & -\gamma_2
  \end{bmatrix},
  \label{eq:LTV_form}
\end{equation}
where
$c_1(t)=-\int_0^1r_1'\bigl(\hat x_2+\vartheta(x_2-\hat x_2)\bigr)d\vartheta$ and
$c_2(t)=\int_0^1r_2'\bigl(\hat x_1+\vartheta(x_1-\hat x_1)\bigr)d\vartheta$.
If $x_2(t)$ and $\hat x_2(t)$ lie in an interval on which
$m_1\le|r_1'|\le M_1$, then $c_1(t)\in[m_1,M_1]$, and similarly for
$c_2$.
\end{lemma}

\begin{proof}
Subtract~\eqref{eq:observer} from~\eqref{eq:GRN} and apply the
fundamental theorem of calculus to
$\vartheta\mapsto r_i(\hat x_j+\vartheta(x_j-\hat x_j))$.
\end{proof}

Lemma~\ref{lem:error} holds for every $C^1$ regulation function; it is
the differential mean-value (LPV) representation used for Lipschitz
systems~\citep{zemouche2013lmi} and is in no way specific to the
logistic family.  What depends on the family is the interval
$[m_i,M_i]$ in which the chord slopes lie: $[0,M_i]$ globally for all
families, and $[m_i,M_i]$ with $m_i>0$ on sets where the slope is
bounded away from zero.

\subsection{Global and set-based certificates}

For $0\le m_i\le M_i$ let
$\mathcal{V}=\{m_1,M_1\}\times\{m_2,M_2\}$ and consider the four
matrix inequalities
\begin{equation}
  \tilde A(v)^TP+P\tilde A(v)+\alpha P\preceq0,\qquad v\in\mathcal{V},
  \label{eq:polytopic_LMI}
\end{equation}
in the unknowns $P\succ0$ and $\alpha>0$.

\begin{theorem}[Convergence and noise attenuation]
\label{thm:cert}
Define
\begin{equation}
  \gamma_{nl}(P,\alpha)=
  \frac{2\sqrt{\mathbf{L}^TP\mathbf{L}}}{\alpha\sqrt{\lambda_{\min}(P)}} .
  \label{eq:gamma_nl}
\end{equation}
\begin{enumerate}
\item[\textup{(i)}] \emph{Global.}  Let Assumption~\ref{ass:reg}
  hold.  If~\eqref{eq:polytopic_LMI} holds with $m_1=m_2=0$ and the
  global bounds $M_1,M_2$ of Assumption~\ref{ass:reg}, then for all
  $\mathbf{x}(0),\hat{\mathbf{x}}(0)\in\mathbb{R}^2$ and all $t\ge0$
  \begin{equation}
    \|\mathbf{e}(t)\|\le\sqrt{\kappa(P)}\,\|\mathbf{e}(0)\|\,e^{-\alpha t/2}
    +\gamma_{nl}(P,\alpha)\,W .
    \label{eq:bound}
  \end{equation}
\item[\textup{(ii)}] \emph{Set-based.}  Let $I_1,I_2$ be closed
  intervals with $\mathbf{x}(t)\in I_1\times I_2$ for all $t\ge0$, let
  $\rho>0$, let $I_i^\rho$ be the $\rho$-neighbourhood of $I_i$, let
  $r_1$ be $C^1$ and monotone on $I_2^\rho$ and $r_2$ on $I_1^\rho$
  (Assumption~\ref{ass:reg} is not needed), and
  let $m_1=\min_{I_2^\rho}|r_1'|$, $M_1=\max_{I_2^\rho}|r_1'|$,
  $m_2=\min_{I_1^\rho}r_2'$, $M_2=\max_{I_1^\rho}r_2'$.
  If~\eqref{eq:polytopic_LMI} holds with these bounds and
  \begin{equation}
    \sqrt{\kappa(P)}\,\|\mathbf{e}(0)\|+\gamma_{nl}(P,\alpha)\,W<\rho ,
    \label{eq:rho_cond}
  \end{equation}
  then~\eqref{eq:bound} holds for all $t\ge0$.
\end{enumerate}
\end{theorem}

\begin{proof}
Let $V=\mathbf{e}^TP\mathbf{e}$.  As long as
$(c_1(t),c_2(t))\in[m_1,M_1]\times[m_2,M_2]$, $\tilde A(c_1,c_2)$ is a
convex combination of the vertex matrices because it is affine in
$(c_1,c_2)$, and~\eqref{eq:polytopic_LMI} with~\eqref{eq:LTV_form}
gives
$\dot V\le-\alpha V+2|w|\,|\mathbf{L}^TP\mathbf{e}|
\le-\alpha V+2W\sqrt{\mathbf{L}^TP\mathbf{L}}\sqrt V$.
Hence $u=\sqrt V$ satisfies
$\dot u\le-(\alpha/2)u+W\sqrt{\mathbf{L}^TP\mathbf{L}}$ where $V>0$; at
points where $V=0$ the same inequality holds for $u_\epsilon=\sqrt{V+\epsilon}$
up to the additional term $(\alpha/2)\sqrt\epsilon$, which vanishes as
$\epsilon\downarrow0$, so by the comparison
lemma
$u(t)\le u(0)e^{-\alpha t/2}+(2W/\alpha)\sqrt{\mathbf{L}^TP\mathbf{L}}$,
and~\eqref{eq:bound} follows from
$\lambda_{\min}(P)\|\mathbf{e}\|^2\le V\le\lambda_{\max}(P)\|\mathbf{e}\|^2$.
In case~(i) the slope condition holds for all times by
Assumption~\ref{ass:reg}, and the plant and observer solutions exist
for all $t\ge0$ because the vector fields are globally Lipschitz.  In case~(ii) it holds as long as
$\|\mathbf{e}(t)\|\le\rho$, because then $\hat x_i(t)\in I_i^\rho$ and
the segment between $x_i$ and $\hat x_i$ lies in $I_i^\rho$.
By~\eqref{eq:rho_cond}, $\|\mathbf{e}(0)\|<\rho$, and on any interval
$[0,t_1]$ on which $\|\mathbf{e}\|\le\rho$ the
bound~\eqref{eq:bound} gives $\|\mathbf{e}\|<\rho$; by continuity
$\|\mathbf{e}(t)\|<\rho$ for all $t\ge0$.
\end{proof}

In part~(i) the regulation functions are understood as extended to
$\mathbb{R}$ (Section~\ref{sec:model}).  For the concentration model
the plant starts in $[0,\infty)^2$ and stays there, while the observer
may be initialised anywhere in $\mathbb{R}^2$ and is not required to
remain non-negative.  Part~(ii) does not assume that the observer
state remains in a prescribed box: the estimate may leave
$\mathcal{B}_0$ during the transient (and even become negative), and
condition~\eqref{eq:rho_cond} accounts for this.  Natural choices for
$I_1\times I_2$ are $\mathcal{B}_0$ or, once the plant has entered it,
any of the forward invariant boxes $\mathcal{B}_k$ of
Lemma~\ref{lem:boxes}.  Part~(i) is a quadratic certificate for
monotone nonlinearities with slopes in $[0,M_i]$, in the spirit of
circle-criterion observers~\citep{arcak2001nonlinear}; it applies to
all four families with the corresponding $M_i$ from
Table~\ref{tab:families} (for Hill exponents $n>1$; see
Section~\ref{sec:subunit} for $n\le1$).  Part~(ii) uses the regulation
functions only on the intervals $I_i^\rho$, whatever their behaviour
elsewhere; this is what allows it to be applied in
Section~\ref{sec:subunit} to Hill functions that are not differentiable
at zero.

\begin{remark}[Relation to existing observer designs]
\label{rem:designs}
Theorem~\ref{thm:cert} is an instance of known constructions, and it
is useful to say which.  Part~(i) is a common quadratic Lyapunov
function for a Lipschitz system written in differential mean-value
form~\citep{rajamani1998observers,zemouche2013lmi}; because the
nonlinearities are monotone it is also covered by the circle-criterion
observers of~\citep{arcak2001nonlinear}, which exploit the sector
$[0,M_i]$ exactly as we do.  Part~(ii) is the corresponding
semi-global statement on a prescribed set.  High-gain observers for
uniformly observable
systems~\citep{gauthier1992simple,khalil2017highgain} take a different
route: the map $\Phi$ of Section~\ref{sec:observability} brings
system~\eqref{eq:GRN} to the triangular form they require, and they
achieve an arbitrarily fast rate, but they need a positive lower bound
on $|r_1'|$ on the whole region of interest, that is, $m_1>0$.  No
family provides such a bound globally, since every sigmoid has a
vanishing slope at saturation, so a high-gain design is semi-global
here, and its correction term in the original coordinates scales like
$1/(\kappa_1m_1)$.  The Luenberger form with fixed gains studied in
this paper gives up arbitrary speed, as Proposition~\ref{prop:limit}
makes precise, in exchange for a global statement with moderate noise
amplification.
\end{remark}

\subsection{Ceiling on the certified rate}

\begin{proposition}[Ceiling on the certified rate]
\label{prop:limit}
Let $T=\gamma_1+\gamma_2+l_1>0$ and
$D(c_1,c_2)=(\gamma_1+l_1)\gamma_2+\kappa_1c_1(\kappa_2c_2-l_2)$.
If~\eqref{eq:polytopic_LMI} holds, then
$\alpha/2\le\bar\eta:=\min_{v\in\mathcal{V}}\eta(v)$ with
\begin{equation}
  \eta(c_1,c_2)=\tfrac T2-\sqrt{\max\{T^2/4-D(c_1,c_2),0\}} .
  \label{eq:eta}
\end{equation}
If $l_2\le\kappa_2m_2$, then $\bar\eta=\eta(m_1,m_2)$.  For $m_1=0$,
$\bar\eta\le\min\{\gamma_2,\gamma_1+l_1\}$, with equality if
$l_2\le\kappa_2m_2$.  If moreover $m_2<M_2$ and $\gamma_1+l_1>\gamma_2$,
then~\eqref{eq:polytopic_LMI} is infeasible at $\alpha=2\gamma_2$, so
the ceiling $\gamma_2$ is not attained.
\end{proposition}

\begin{proof}
Inequality~\eqref{eq:polytopic_LMI} at a vertex $v$ implies that the
eigenvalues of $\tilde A(v)$ have real part at most $-\alpha/2$, and
$-\eta(v)$ is the largest real part of the roots of
$\xi^2+T\xi+D(v)$.  Since $\eta$ is non-decreasing in $D$, and $D$ is
non-decreasing in $c_1$ and $c_2$ when $\kappa_2c_2-l_2\ge0$, the
minimum is at $(m_1,m_2)$.  For $c_1=0$, $\tilde A$ is triangular with
eigenvalues $-(\gamma_1+l_1)$ and $-\gamma_2$.  At $c_1=0$ and
$\alpha=2\gamma_2$ the $(2,2)$ entry of the left-hand side
of~\eqref{eq:polytopic_LMI} vanishes, so negative semidefiniteness
forces its $(1,2)$ entry,
$(\gamma_2-\gamma_1-l_1)p_{12}+(\kappa_2c_2-l_2)p_{22}$, to vanish as
well, i.e.\ $p_{12}/p_{22}=(\kappa_2c_2-l_2)/(\gamma_1+l_1-\gamma_2)$.
The vertices $c_2=m_2$ and $c_2=M_2$ require different values of this
ratio, so no common $P$ exists.
\end{proof}

Proposition~\ref{prop:limit} states precisely what a positive lower
slope bound buys.  With $m_1=0$ the unmeasured species is not
observable at the corresponding vertex, only detectable, and the rate
of a certificate of the form~\eqref{eq:polytopic_LMI} is at most
$\min\{\gamma_2,\gamma_1+l_1\}$: whatever the gains, it cannot exceed
the degradation rate $\gamma_2$, the ceiling equals $\gamma_2$ only if
$\gamma_1+l_1\ge\gamma_2$, as for the baseline gains, and even then it
is approached but not attained.  How closely it is approached depends
on the family and the gains (Section~\ref{sec:compare}).  This is the
situation of \emph{every} family in the global case, the logistic
included.  A certified rate above $\gamma_2$ requires $m_1>0$,
and~\eqref{eq:eta} quantifies the gain; it is modest when $m_1$ is
small.

\subsection{Computing and certifying the rate}

The supremal rate for which~\eqref{eq:polytopic_LMI} is feasible is
\begin{equation*}
  \alpha^{\star}=-\inf_{P\succ0}\,\max_{v\in\mathcal{V}}\,
  \lambda_{\max}\bigl(P^{-1/2}(\tilde A(v)^TP+P\tilde A(v))P^{-1/2}\bigr),
\end{equation*}
a quasi-convex generalised eigenvalue problem~\citep{boyd1994linear}.
We do not rely on a numerical optimiser to locate $\alpha^\star$: every
reported rate is bracketed by two certificates that can be checked by
evaluating $2\times2$ eigenvalue problems.

\emph{Lower bounds.}  A rate $\alpha$ is certified when a matrix
$P\succ0$, printed with finitely many digits, satisfies
$\max_{v}\lambda_{\max}(\tilde A(v)^TP+P\tilde A(v)+\alpha P)<0$.
Reported rates, written $\alpha_{\mathrm{cert}}/2$, are rounded down
so that the printed $P$ passes this strict check.

\emph{Upper bounds.}  Besides Proposition~\ref{prop:limit}, we use the
theorem of alternatives for linear matrix
inequalities~\citep{boyd1994linear} in the following form.

\begin{lemma}[Dual bound on the rate]
\label{lem:dual}
Let $\alpha\ge0$.  If there are matrices $Z_v\succeq0$,
$v\in\mathcal{V}$, not all zero, such that
\begin{equation}
  G=\sum_{v\in\mathcal{V}}\bigl(\tilde A(v)Z_v+Z_v\tilde A(v)^T
  +\alpha Z_v\bigr)\succeq0,
  \label{eq:dual}
\end{equation}
then $\alpha^\star\le\alpha$.
\end{lemma}

\begin{proof}
Suppose some $P\succ0$ satisfies $\tilde A(v)^TP+P\tilde A(v)+\alpha
P\prec0$ for every $v$.  Then
$\sum_v\textup{tr}\bigl((\tilde A(v)^TP+P\tilde A(v)+\alpha P)Z_v\bigr)<0$,
whereas this sum equals $\textup{tr}(PG)\ge0$; so no such $P$ exists.  If
some $P\succ0$ satisfied~\eqref{eq:polytopic_LMI} at a rate
$\alpha'>\alpha$, then
$\tilde A(v)^TP+P\tilde A(v)+\alpha P\preceq-(\alpha'-\alpha)P\prec0$
for every $v$, which has just been excluded.  Hence
$\alpha^\star\le\alpha$.
\end{proof}

With $\alpha=0$, Lemma~\ref{lem:dual} proves that the four vertex
matrices have no common quadratic Lyapunov function at all.

\emph{Procedure.}  Candidates are found numerically and then checked.
The objective defining $\alpha^\star$ is unchanged when $P$ is
multiplied by a positive number, so the search uses the two-parameter
representative
$\bigl(\begin{smallmatrix}e^{u}&w\\w&e^{-u}\end{smallmatrix}\bigr)$,
$|w|<1$ (determinant $1-w^2$), with a grid search followed by
Nelder--Mead refinement (R function \texttt{optim}); this
parametrisation keeps the problem well scaled when the gains are large,
and accepted matrices are rescaled to $\det P=1$.  Candidate dual
matrices are obtained from the primal optimum by complementary
slackness: they are rank one, along the null vectors of the active
residuals.  The test~\eqref{eq:dual} is invariant under scaling of the
$Z_v$; we normalise $\sum_v\textup{tr}\,Z_v=1$ and write $Z_v=L_vL_v^T$
with $L_v$ lower triangular, so that $Z_v\succeq0$ holds exactly for any
printed entries, and $\lambda_{\min}(G)>0$ is checked with the printed
(rounded) factors.  Every rate in Tables~\ref{tab:compare}
and~\ref{tab:subunit} is bracketed in this way within
$10^{-4}$~min$^{-1}$, and the certificates are listed in
\ref{app:certificates}.

\subsection{Baseline certificates}

With the gains above and $M_1=M_2=\lambda/4$, the global
inequality~\eqref{eq:polytopic_LMI} ($m_i=0$) is feasible with
$\alpha_{\mathrm{cert}}/2=0.0999$~min$^{-1}$, within $0.1\%$ of the
ceiling $\gamma_2=0.1$~min$^{-1}$ of Proposition~\ref{prop:limit},
which is itself infeasible; this near-optimal certificate is listed in
\ref{app:certificates}.  For the trajectory bound we use instead a
simpler and more conservative certificate with round entries,
\begin{equation}
  P=\begin{bmatrix}11.4&2.7\\2.7&1.7\end{bmatrix},\qquad
  \alpha=0.18~\textup{min}^{-1},
  \label{eq:P_cert}
\end{equation}
for which the largest eigenvalue of the left-hand side
of~\eqref{eq:polytopic_LMI} over the four vertices is $-0.019$,
$\kappa(P)\approx12.1$, so that $\sqrt{\kappa(P)}\approx3.48$, and
$\gamma_{nl}\approx53.4$.  Hence, for all initial conditions,
$\|\mathbf{e}(t)\|\le3.48\,\|\mathbf{e}(0)\|e^{-0.09t}+53.4\,W$, with
the rate $\alpha/2=0.09$~min$^{-1}$ of this particular certificate.

For the set-based version with $I_i=[0,10]$ (the invariant box
$\mathcal{B}_0$), $m_i=\lambda g(-\rho)$, and $\alpha^\star/2$ agrees
with the ceiling $\bar\eta$ to within $10^{-4}$~min$^{-1}$: it is
$0.111$~min$^{-1}$ for $\rho\to0$, $0.104$ for $\rho=1$~nM and $0.100$
for $\rho=5$~nM.  The certified improvement over the global rate is
thus at most about $11\%$ for this parameter set, because $\mu$ is
small.  It also comes with a restricted basin: with the matrices $P$ of
Table~\ref{tab:certs}, condition~\eqref{eq:rho_cond} with $W=0$ admits
$\|\mathbf{e}(0)\|<0.16$~nM for $\rho=1$ and
$\|\mathbf{e}(0)\|<0.71$~nM for $\rho=5$.

\begin{remark}[Local rate and linearised noise gain]
\label{rem:local}
The certificates above are worst-case over the slope polytope and are
conservative near the equilibrium, where the error dynamics are
governed by $A_{cl}$: by Lyapunov's indirect
method~\citep{khalil2002nonlinear} the error converges locally at any
rate below the placed poles ($1$~min$^{-1}$ here).  Likewise, for the
linearisation $\dot{\mathbf{e}}=A_{cl}\mathbf{e}-\mathbf{L}w$ the
steady-state error is bounded by $\gamma_{pp}W$ with the peak-to-peak
bound $\gamma_{pp}=\int_0^\infty\|e^{A_{cl}t}\mathbf{L}\|dt$; for the
baseline $A_{cl}+I$ is nilpotent,
$e^{A_{cl}t}=e^{-t}(I+(A_{cl}+I)t)$, and $\gamma_{pp}\approx2.70$.
The global gain $\gamma_{nl}\approx53.4$ is therefore conservative by
more than an order of magnitude for small errors; it is the price of a
statement that holds for all initial conditions and all admissible
slopes.  Faster placed poles increase $|l_2|$ quadratically and the
noise gains with it, the usual trade-off of high-gain
designs~\citep{khalil2017highgain}.
\end{remark}

\section{Numerical illustration}
\label{sec:numerics}

\begin{figure*}[t]
  \centering
  \includegraphics[width=0.85\linewidth]{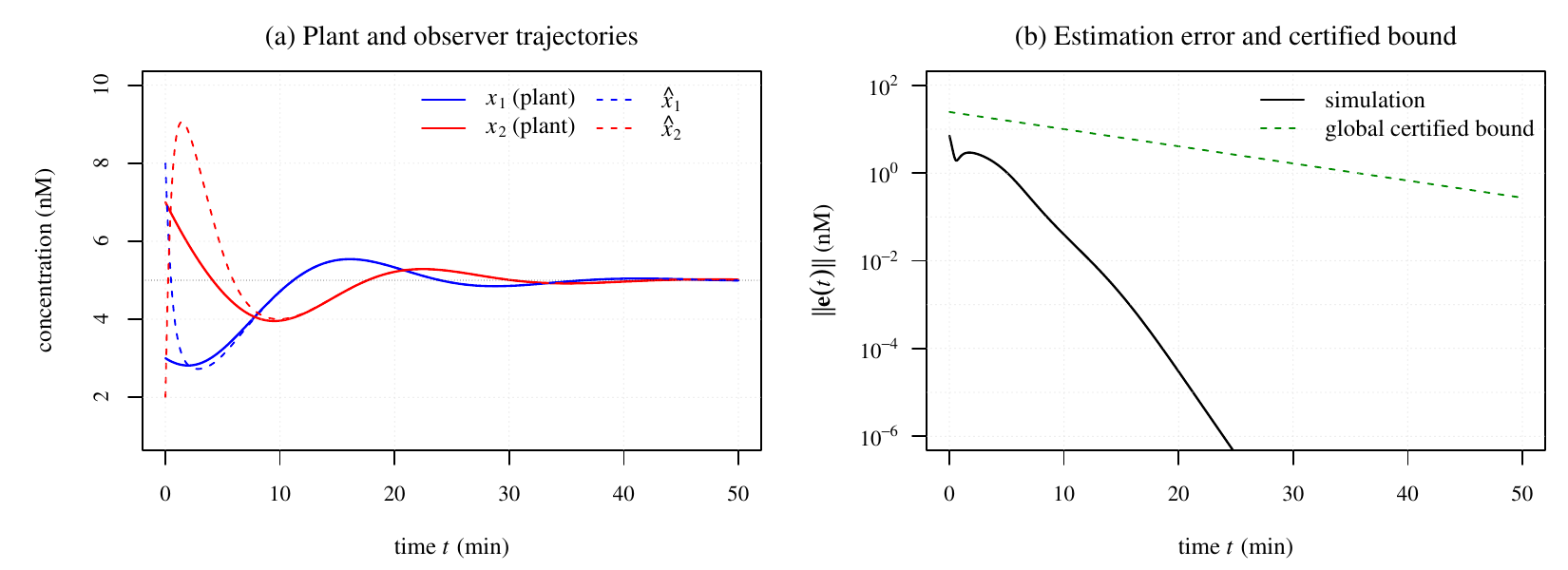}
  \caption{Baseline observer without noise.  (a) Plant (solid) and
    observer (dashed) trajectories; the dotted line marks the
    equilibrium value $5$~nM.  (b) Estimation error (solid) and the
    global certified bound $3.48\,\|\mathbf{e}(0)\|e^{-0.09t}$ of
    Theorem~\ref{thm:cert}(i) with~\eqref{eq:P_cert} (dashed).}
  \label{fig:obs}
\end{figure*}

All computations were done in R~4.3.3 with the package
\texttt{deSolve}~1.40; noise-free trajectories were integrated with
\texttt{lsoda} (relative tolerance $10^{-11}$, absolute tolerance
$10^{-13}$).  We
integrate~\eqref{eq:GRN} and~\eqref{eq:observer} for the baseline
with $\mathbf{x}(0)=(3,7)$ and $\hat{\mathbf{x}}(0)=(8,2)$~nM
($\|\mathbf{e}(0)\|=\sqrt{50}\approx7.07$~nM) and $w=0$
(Fig.~\ref{fig:obs}).  The error falls below $1$, $0.5$, $0.1$ and
$0.01$~nM at $t\approx5.05$, $6.16$, $8.51$ and $12.29$~min.  Its
logarithmic decay rate grows from about $0.6$~min$^{-1}$ for
$t\in[8,16]$ to about $0.9$~min$^{-1}$ for $t\in[24,32]$, approaching
the placed-pole rate as the plant settles near the equilibrium.  The
global bound of Theorem~\ref{thm:cert}(i) with~\eqref{eq:P_cert} is
respected throughout and is, as expected, conservative.
As Lemma~\ref{lem:error} requires, integrating~\eqref{eq:LTV_form}
with the chord slopes computed along the trajectory reproduces the
nonlinear error; the discrepancy, below $10^{-7}$~nM, comes from the
interpolation of the stored trajectory.

With measurement noise we start from $\mathbf{x}(0)=(5,5)$,
$\hat{\mathbf{x}}(0)=(5.1,4.9)$~nM.  The noise is piecewise constant,
and the integration (\texttt{lsoda}, relative tolerance $10^{-9}$,
absolute tolerance $10^{-11}$) is restarted at every jump, so that the
solver never steps across a discontinuity; the error is sampled every
$0.05$~min.  Random noise is uniform on $[-W,W]$, drawn with R's
default Mersenne-Twister generator and held constant on each interval
$[0.1k,0.1(k+1))$~min ($20$ realisations with seeds $1$ to $20$,
horizon $100$~min, peak of $\|\mathbf{e}\|$ over $t\ge30$~min); the
worst sampled peak is $1.17$--$1.18\,W$ for $W\in\{0.05,0.2,1\}$~nM.
Square-wave noise equal to $+W$ on the first half period and
alternating in sign, with periods $4,6,8,10,14,20$ and $30$~min
(horizon $200$~min, peak over $t\ge60$~min), gives at most
$2.44$--$2.49\,W$.  Both are below the linearised bound
$\gamma_{pp}W=2.70\,W$ of Remark~\ref{rem:local} and far below the
global bound $53.4\,W$.  These are the worst of the sampled inputs, not worst-case values;
they illustrate the bounds and are not a substitute for them.

\section{Comparison of regulation families}
\label{sec:compare}

\subsection{Matched parameter sets}

We compare the four families of Table~\ref{tab:families} (Hill
exponents $n>1$) with
$\kappa_i=1$, $\gamma_i=0.1$ (so that $\mathcal{B}_0=[0,10]^2$), the
same family for $r_1$ and $r_2$, and parameters matched so that the
sigmoidal part of every activation function is half-maximal at the same
point, with the same slope there.  For the Hill function,
$h(K)=1/2$ and $h'(K)=n/(4K)$, while $f(\theta)=1/2$ and
$f'(\theta)=\lambda/4$, so the matching is
\begin{equation}
  \theta=K,\qquad \lambda=\frac{n}{K};
  \label{eq:matching}
\end{equation}
for the shifted Hill function, where $x_0$ is the input shift and $K$
the half-response point in the shifted coordinate $x+x_0$, it reads
$\theta=K-x_0$, $\lambda=n/K$.  In the leaky case the parameters of
the underlying Hill sigmoid are matched in the same way, but the
complete response then has value $(1+b)/2$ and slope $(1-b)n/(4K)$ at
$K$, so neither is matched exactly.  This is a two-condition local
calibration: it fixes the threshold and the slope there, and nothing
else.  The parameter sets are:
\begin{itemize}
\item \emph{Scenario A} (threshold in mid-range, slope $1/4$):
  logistic $\theta=5$, $\lambda=1$; Hill $K=5$, $n=5$; leaky Hill
  $K=5$, $n=5$, $b=0.05$; shifted Hill $x_0=1$, $K=6$, $n=6$.
\item \emph{Scenario B} (threshold low in the range, slope $1/2$):
  logistic $\theta=2$, $\lambda=2$; Hill $K=2$, $n=4$; leaky Hill
  $K=2$, $n=4$, $b=0.05$; shifted Hill $x_0=1$, $K=3$, $n=6$.
\end{itemize}

\subsection{Computed quantities}

\begin{figure*}[t]
  \centering
  \includegraphics[width=0.85\linewidth]{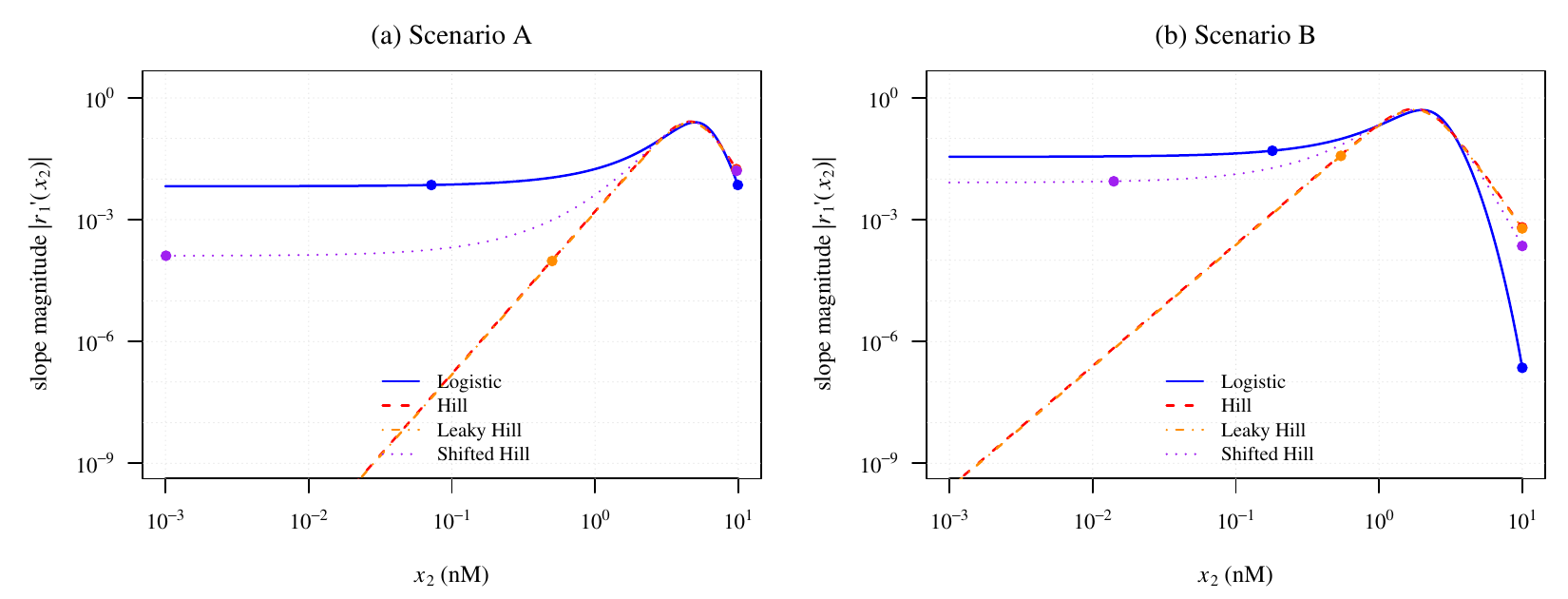}
  \caption{Slope magnitude $|r_1'(x_2)|$ of the four families in
    scenarios A and B.  Dots mark the ends of the operating interval
    $J_2$ of each model (where they fall inside the plotted range).
    The Hill and leaky-Hill curves differ only by the factor $1-b$;
    their operating intervals differ.}
  \label{fig:slopes}
\end{figure*}

For each model we compute the operating set of Lemma~\ref{lem:boxes},
the lower slope bounds $m_1$ on $J_2$ and $m_2$ on $J_1$, the gains of
Proposition~\ref{prop:gains} placing both poles at $-1$ at the model's
own equilibrium, the ceiling $\bar\eta$ of Proposition~\ref{prop:limit}
for the slope ranges of the operating set, and the certified rate
$\alpha_{\mathrm{cert}}/2$ of~\eqref{eq:polytopic_LMI} on the vertices
of these ranges, together with that of the global certificate of
Theorem~\ref{thm:cert}(i) (Table~\ref{tab:compare},
Fig.~\ref{fig:slopes}).  The condition $l_2\le\kappa_2m_2$ holds in all
rows, so $\bar\eta=\eta(m_1,m_2)$.  By Theorem~\ref{thm:obs} with $\tau=1$~min,
$\sigma_{\min}\ge s(\tau\kappa_1m_1)$, and
$s(\tau\kappa_1m_1)=0.995\,\tau\kappa_1m_1$ to three digits in all
rows, the factor being $1/\sqrt{1+\tau^2\gamma_1^2}$.

The minima $m_1,m_2$ on $J_2,J_1$ are operating-set conditioning
indicators, and the set-based rates are those of
Theorem~\ref{thm:cert}(ii) in the limit $\rho\to0$.  Two qualifications
apply.  First, $\mathcal{B}_\infty$ is attracting but is reached only
asymptotically, so a set-based certificate is used in practice with a
forward invariant box $\mathcal{B}_k$ slightly larger than
$\mathcal{B}_\infty$.  Second, a rate obtained with $\rho=0$ is a limiting
benchmark, not itself a guarantee: Theorem~\ref{thm:cert}(ii) requires
$\rho>0$, with the slope bounds taken on the enlarged intervals, and
condition~\eqref{eq:rho_cond}.  The enlarged polytope contains the
$\rho=0$ one, so the finite-$\rho$ rate is no larger, and the minima
vanish as soon as the enlarged intervals reach the constant extension
of a Hill-type function.  For the baseline, the corresponding rates
and admissible initial errors are given in Section~\ref{sec:design}.

\begin{table*}[t]
\centering
\caption{Operating intervals $J_1,J_2$ (nM), lower slope bounds
$m_2=\min_{J_1}r_2'$ and $m_1=\min_{J_2}|r_1'|$ with the end of $J_2$
where $m_1$ is attained (slopes in nM$^{-1}$), gain $l_2$ for both
poles at $-1$ ($l_1=1.8$ in all rows), rate ceiling $\bar\eta$ of
Proposition~\ref{prop:limit} on the operating set, and certified
rates $\alpha_{\mathrm{cert}}/2$ on the operating set (``set'') and
globally (``global'', $m_i=0$), in min$^{-1}$.  The set-based rates
are evaluated with $\rho=0$, the limit of Theorem~\ref{thm:cert}(ii)
as $\rho\downarrow0$ (see text).  Certified rates are rounded down;
each lies within $10^{-4}$ below $\alpha^\star/2$
(\ref{app:certificates}).}
\label{tab:compare}
\small
\begin{tabular}{@{}llllllrccc@{}}
\toprule
 & & & & & & & & \multicolumn{2}{c}{$\alpha_{\mathrm{cert}}/2$} \\
\cmidrule(l){9-10}
Family & $J_1$ & $m_2$ & $J_2$ & $m_1$ & at & $l_2$ & $\bar\eta$ & set & global \\
\midrule
\multicolumn{10}{@{}l}{\emph{Scenario A}}\\
Logistic     & $[0.072,\,9.93]$             & $7.1{\times}10^{-3}$ & $[0.072,\,9.93]$             & $7.1{\times}10^{-3}$  & both  & $-2.990$ & $0.1120$ & $0.1119$ & $0.0999$ \\
Hill         & $[0.35,\,10.0]$              & $2.4{\times}10^{-5}$ & $[1.7{\times}10^{-5},\,9.70]$ & $1.4{\times}10^{-22}$ & lower & $-2.990$ & $0.1000$ & $0.0999$ & $0.0999$ \\
Leaky Hill   & $[0.83,\,10.0]$              & $7.3{\times}10^{-4}$ & $[0.50,\,9.71]$              & $9.6{\times}10^{-5}$  & lower & $-3.272$ & $0.1002$ & $0.1001$ & $0.0999$ \\
Shifted Hill & $[0.29,\,10.0]$              & $4.7{\times}10^{-4}$ & $[1.0{\times}10^{-3},\,9.74]$ & $1.3{\times}10^{-4}$  & lower & $-2.990$ & $0.1002$ & $0.1001$ & $0.0999$ \\
\midrule
\multicolumn{10}{@{}l}{\emph{Scenario B}}\\
Logistic     & $[1.1{\times}10^{-6},\,9.74]$ & $3.8{\times}10^{-7}$ & $[0.18,\,10.0]$              & $2.3{\times}10^{-7}$  & upper & $-2.700$ & $0.1000$ & $0.0994$ & $0.0994$ \\
Hill         & $[0.016,\,10.0]$             & $1.0{\times}10^{-6}$ & $[4.2{\times}10^{-8},\,9.98]$ & $1.8{\times}10^{-23}$ & lower & $-3.976$ & $0.1000$ & $0.0994$ & $0.0994$ \\
Leaky Hill   & $[0.52,\,9.95]$              & $6.2{\times}10^{-4}$ & $[0.54,\,9.98]$              & $6.1{\times}10^{-4}$  & upper & $-6.350$ & $0.1022$ & $0.0999$ & $0.0978$ \\
Shifted Hill & $[0.0041,\,9.99]$            & $2.3{\times}10^{-4}$ & $[0.014,\,10.0]$             & $2.2{\times}10^{-4}$  & upper & $-3.501$ & $0.1004$ & $0.0999$ & $0.0995$ \\
\bottomrule
\end{tabular}
\end{table*}

\subsection{Observations}

\emph{(1) The unshifted Hill model is observable on its operating set,
but with a negligible margin.}  The operating interval excludes
$x_2=0$, so $m_1>0$, as the rank condition requires; but the lower end
of $J_2$ is of order $10^{-5}$ to $10^{-8}$~nM and $m_1$ of order
$10^{-22}$ to $10^{-23}$.  The statement that the Hill observability
matrix ``loses rank at basal expression'' is exact only on the line
$x_2=0$, which solutions starting elsewhere do not reach; the
practically relevant statement is the size of $m_1$.

\emph{(2) Basal expression or an input shift removes most of that
degeneracy.}  With a $5\%$ leak or a shift of $1$~nM, $m_1$ is of order
$10^{-4}$, eighteen to nineteen orders of magnitude larger than for
the unshifted Hill model.  The comparison that matters for the
logistic function is therefore with these generalised Hill forms.

\emph{(3) Neither the logistic nor the generalised Hill forms
dominate.}  In scenario A the logistic bound exceeds those of the leaky
and shifted Hill models by factors of about $75$ and $55$.  In
scenario B the ranking is reversed: the logistic slope at the upper end
of $J_2$ is $2.3\times10^{-7}$, about three orders of magnitude below
the generalised Hill values.  The reason is the tails.  The logistic
slope is nearly constant, equal to about $\lambda e^{-\lambda\theta}$,
at the basal end but decays like $e^{-\lambda(x-\theta)}$ above the
threshold, while Hill-type slopes vanish (or become small) at the
basal end and decay like $x^{-(n+1)}$ above the threshold.  Which model
is better conditioned depends on where the threshold sits within the
range $[0,\kappa_2/\gamma_2]$ of the unmeasured species: roughly, the
logistic form is favoured when the range extends only a few multiples
of $1/\lambda$ above the threshold, and Hill-type forms otherwise.

\emph{(4) The differences concern conditioning much more than
certified convergence rates.}  In all eight cases the certified rate
$\alpha_{\mathrm{cert}}/2$ on the operating set lies between $0.0994$ and
$0.1119$~min$^{-1}$, within $0.003$~min$^{-1}$ of the ceiling
$\bar\eta$.  The global certificate of Theorem~\ref{thm:cert}(i) is
feasible for every model, with certified rates between $0.0994$ and
$0.0999$~min$^{-1}$, close to but strictly below the global ceiling
$\gamma_2=0.1$, except for the leaky Hill model of scenario~B, whose
larger gain $|l_2|$ lowers it to $0.0978$~min$^{-1}$.  What $m_1$ controls is the bound
$1/s(\tau\kappa_1m_1)$ of~\eqref{eq:inverse_lip} on the sensitivity of the
reconstructed state, and the size of the gain $l_2$ when the
equilibrium lies in a low-slope region
(Proposition~\ref{prop:gains}).

\section{Hill exponents not exceeding one}
\label{sec:subunit}

The preceding sections assume $n>1$ for Hill-type functions.  We now
state what changes for $n\le1$, where the roles of the lower and upper
slope bounds are exchanged.

\subsection{Regularity}

For $h(x)=x^n/(K^n+x^n)$ and $x>0$,
$h'(x)=nK^nx^{n-1}/(K^n+x^n)^2$ and
\begin{equation}
  \frac{d}{dx}\ln h'(x)=\frac{(n-1)K^n-(n+1)x^n}{x\,(K^n+x^n)} .
  \label{eq:logderiv}
\end{equation}

\emph{Case $n=1$.}  $h'(x)=K/(K+x)^2$ decreases from $1/K$ at $x=0$.
The slope is bounded and positive at the basal end, and on $[0,B]$ its
minimum is $K/(K+B)^2$.  The constant continuation to $x<0$ is not
$C^1$, but the tangent-line continuation $h(x)=x/K$ for $x<0$ is, with
slope $1/K$; with it Assumption~\ref{ass:reg} holds with $M=1/K$ and
all results of Sections~\ref{sec:observability}--\ref{sec:design}
apply.

\emph{Case $0<n<1$.}  By~\eqref{eq:logderiv}, $h'$ is strictly
decreasing on $(0,\infty)$, with $h'(x)\sim nK^{-n}x^{n-1}\to\infty$ as
$x\to0^+$.  Hence, for every $B>0$,
\begin{equation}
  \inf_{0<x\le B}h'(x)=h'(B)>0,\qquad \sup_{0<x\le B}h'(x)=\infty .
  \label{eq:subunit_bounds}
\end{equation}
The function is H\"older continuous,
$|h(x)-h(y)|\le K^{-n}|x-y|^n$ (from $|a^n-b^n|\le|a-b|^n$ and the fact
that $u\mapsto u/(1+u)$ is $1$-Lipschitz), but it has no $C^1$
continuation across $x=0$ and no finite Lipschitz constant on an
interval containing zero.  A leaky form $b+(1-b)h$ has the same
singularity: the leak changes the basal value, not the slope.  A
shifted form $h(x+x_0)$, $x_0>0$, does not: on $x\ge0$ its slope
decreases from $h'(x_0)<\infty$, and the tangent-line continuation to
$x<0$ gives Assumption~\ref{ass:reg} with $M=h'(x_0)$.

For the plant itself the singularity is harmless.  The construction of
Lemma~\ref{lem:boxes} and the existence and uniqueness of the
equilibrium in Proposition~\ref{prop:equil} use only continuity and
monotonicity of $r_i$.  Solutions starting in
$[0,\infty)^2\setminus\{0\}$ are unique.  Away from the coordinate axes
the vector field is locally Lipschitz.  Near a point of the axis
$x_2=0$ with $x_1>0$, the second component
$\kappa_2r_2(x_1)-\gamma_2x_2$ is positive and locally Lipschitz, so
$x_2$ can serve as independent variable: $dx_1/dx_2$ is then locally
Lipschitz in the dependent variable $x_1$ and continuous in $x_2$,
which gives uniqueness of the orbit, and $t$ is recovered from
$dt/dx_2=1/\dot x_2$.  The axis $x_1=0$ with $x_2>0$ is treated in the
same way with $x_1$ as independent variable.  Since $\dot x_1>0$ near
the origin, no solution starting elsewhere reaches it, and no result
below depends on solutions starting at the origin.

\subsection{Observability}

Let $r_1$ be an unshifted (possibly leaky) Hill repression with
$0<n<1$.

\emph{Injectivity.}  Theorem~\ref{thm:obs}(i) uses strict monotonicity
only, so the observability map is injective on
$\mathbb{R}\times[0,\infty)$.

\emph{Rank condition.}  For $x_2>0$ the slope is finite and non-zero and
the rank condition holds.  At $x_2=0$ the matrix
$\mathcal{O}$ is not defined; this is a failure of differentiability,
not a loss of rank.

\emph{Conditioning.}  The function $s$ of~\eqref{eq:s_def} satisfies
$s(d)\to1$ as $d\to\infty$.  By~\eqref{eq:subunit_bounds}, for
$0<x_2\le B$,
$\sigma_{\min}(\mathcal{O})\ge s(\tau\kappa_1|r_1'(B)|)>0$, and
$\sigma_{\min}\to1$ as $x_2\to0^+$ while $\sigma_{\max}\to\infty$.  The
smallest singular value does not deteriorate at the basal end; it is
the largest one, and with it the condition number, that is unbounded.
The inverse bound survives the loss of differentiability.

\begin{corollary}[Inverse bound without differentiability]
\label{cor:ac}
Let $r_1$ be strictly monotone and absolutely continuous on
$J=[0,B]$, with $|r_1'|\ge m_1>0$ almost everywhere on $J$.
Then~\eqref{eq:inverse_lip} holds for all
$\mathbf{x},\mathbf{x}'\in\mathbb{R}\times J$.  In particular it holds
for an unshifted or leaky Hill repression with $0<n<1$, with
$m_1=|r_1'(B)|$.
\end{corollary}

\begin{proof}
The proof of Theorem~\ref{thm:obs}(iii) uses $r_1'$ only through the
secant slope $\bar c=(r_1(x_2)-r_1(x_2'))/(x_2-x_2')$.  For an
absolutely continuous monotone $r_1$ this is the average of $r_1'$
over the segment, so $|\bar c|\ge m_1$.  The Hill case follows
from~\eqref{eq:subunit_bounds}.
\end{proof}

In the sense of this inverse Lipschitz bound, a Hill function with
$n<1$ is better conditioned near zero than the logistic function; the
condition number of $\mathcal{O}$, in contrast, diverges there.

\subsection{Observer certificates}

The error representation of Lemma~\ref{lem:error} remains valid as an
integral identity while both arguments stay in $[0,\infty)$, because
$h$ is absolutely continuous on compact intervals, but
by~\eqref{eq:subunit_bounds} the chord slopes are not bounded above
when both arguments approach zero.  Consequently:

\emph{Global certificate.}  Theorem~\ref{thm:cert}(i) requires finite
$M_i$ and does not apply to unshifted or leaky Hill functions with
$n<1$.  In addition, since no $C^1$ continuation to negative arguments
exists, the observer~\eqref{eq:observer} is then not covered by the
standard existence and uniqueness theory once its state reaches zero; a
positivity-preserving observer~\citep{proskurnikov2019positive}, a
regularisation of $h$ near zero, or a separate analysis would be
needed.  The theorem does apply to the
logistic function, to Hill functions with $n=1$, and to shifted Hill
functions with any $n>0$.

\emph{Set-based certificate.}  Theorem~\ref{thm:cert}(ii) applies to
unshifted Hill functions with $n<1$ provided the enlarged intervals
satisfy $I_i^\rho\subset(0,\infty)$, with $M_i$ the (finite) slope at
their lower end.  The certificate deteriorates as that end approaches
zero, because the slope range widens.

\subsection{Scenario C: exponent one half}

We take $\kappa_i=1$, $\gamma_i=0.1$ and match at the threshold
by~\eqref{eq:matching}: logistic $\theta=5$, $\lambda=0.1$; Hill $K=5$,
$n=0.5$; leaky Hill $K=5$, $n=0.5$, $b=0.05$; shifted Hill $x_0=1$,
$K=6$, $n=0.6$.  The slope at the threshold is $0.025$ for the
logistic, Hill and shifted Hill models and $(1-b)\times0.025$ for the
leaky one.  The matched curves differ considerably away from the
threshold: the logistic function is almost linear on $[0,10]$, rising
from $0.38$ to $0.62$, whereas the Hill function rises from $0$ to
$0.59$.

With such a shallow regulation the iteration of
Lemma~\ref{lem:boxes} contracts, and for each model the operating set
$\mathcal{B}_\infty$ reduces to that model's own equilibrium: $(5,5)$~nM
for the logistic, Hill and shifted Hill models, and
$(5.18,5.29)$~nM for the leaky model, whose two responses are not
complements.  The operating sets therefore carry no slope information,
and the comparison has to be made on the transient set
$\mathcal{B}_0=[0,10]^2$, which also has to accommodate the excursions
of the observer state.  Table~\ref{tab:subunit} gives the slope bounds on
$[0,10]$ and the certified rates $\alpha_{\mathrm{cert}}/2$ of the
four-vertex inequality~\eqref{eq:polytopic_LMI}, with gains placing
both poles at $-1$ ($l_1=1.8$; $l_2=-32.4$, and $-36.1$ for the leaky
model; the large $|l_2|$ reflects the small slope at the equilibrium).

\begin{table}[ht]
\centering
\caption{Scenario C ($n<1$).  Slope bounds of $r_1$ on $[0,10]$, and
certified rates $\alpha_{\mathrm{cert}}/2$ (min$^{-1}$) of the global
certificate (Theorem~\ref{thm:cert}(i)) and of the set-based
certificate (Theorem~\ref{thm:cert}(ii) with $\rho=0$, the limit as
$\rho\downarrow0$) on the interval shown.  Slopes in nM$^{-1}$;
certified rates $\alpha_{\mathrm{cert}}/2$, rounded down, each within
$10^{-4}$ below $\alpha^\star/2$ (\ref{app:certificates}).}
\label{tab:subunit}
\small
\begin{tabular}{@{}lllll@{}}
\toprule
Family & $\min|r_1'|$ & $\sup|r_1'|$ & Global & Set-based \\
\midrule
Logistic     & $0.0235$ & $0.0250$ & $0.0999$ & $0.7794$ on $[0,10]$ \\
Hill         & $0.0121$ & $\infty$ & n/a      & $0.2741$ on $[1,10]$ \\
Leaky Hill   & $0.0115$ & $\infty$ & n/a      & $0.2676$ on $[1,10]$ \\
Shifted Hill & $0.0132$ & $0.114$  & $0.0929$ & $0.2698$ on $[0,10]$ \\
\bottomrule
\end{tabular}
\end{table}

For the Hill and leaky-Hill models the set-based inequality is feasible
on $[1,10]$, where the slopes lie in $[0.0121,0.107]$ and $[0.0115,0.101]$
respectively.  For these gains it is infeasible on $[0.1,10]$, where
the upper slope is $0.54$ and $0.52$, although all four vertex
matrices are Hurwitz.  On intervals $[\varepsilon,10]$ the slope
polytope grows as $\varepsilon$ decreases, so feasibility is monotone
in $\varepsilon$, and the transition is located by certificates on
both sides (\ref{app:certificates}): the inequality is feasible for
$\varepsilon\ge0.368$ and, by a dual certificate at $\alpha=0$,
infeasible for $\varepsilon\le0.367$ for the Hill model, so the
transition lies in $(0.367,0.368)$; for the leaky model it lies in
$(0.391,0.392)$.  The upper
slope at the transition is about $0.23$ and $0.21$ respectively.  The square $[1,10]^2$ is forward
invariant for these two plants ($\dot x_1=r_1(x_2)-0.1x_1>0$ on
$x_1=1$, since $r_1\ge0.41$ on $[1,10]$, and likewise $\dot x_2>0$ on
$x_2=1$), so the set-based certificate on $[1,10]$ applies to plant
trajectories that start in this square, or from the time they enter
it, with condition~\eqref{eq:rho_cond} checked at that time.  The logistic slopes vary by $6\%$
over the whole box, and its set-based rate $0.7794$ equals the ceiling
$\bar\eta=0.7795$ of Proposition~\ref{prop:limit} to rounding, not far
below the placed pole.

These set-based rates are $\rho\downarrow0$ benchmarks, and the finite-$\rho$
version shows what they cost.  With $\rho=0.5$~nM (slopes taken on
$[-0.5,10.5]$, or on $[0.5,10.5]$ for the Hill and leaky models), the
certified rates fall to $0.7584$ (logistic), $0.0911$ (Hill), $0.0738$
(leaky Hill) and $0.2621$ (shifted Hill)~min$^{-1}$, and the matrices
$P$ that achieve them are so ill-conditioned that
condition~\eqref{eq:rho_cond} with $W=0$ admits only
$\|\mathbf{e}(0)\|<0.003$, $0.023$, $0.021$ and $0.017$~nM
(Table~\ref{tab:frho}).  The fast
set-based rates are therefore local statements; the global
certificate, with rate close to $\gamma_2$, is the one that applies to
arbitrary initial errors.

Three conclusions follow for $n<1$.  First, the comparison of lower
slope bounds made in Section~\ref{sec:compare} does not carry over: near
zero the unshifted Hill slope is large, not small, and the inverse of
its observability map has a small Lipschitz constant there.  Second,
what the logistic function offers in this regime is a finite,
threshold-independent upper slope bound on the whole real line, hence
a global certificate that includes the basal boundary and an observer
that is well defined wherever its state goes; unshifted Hill
certificates are confined to regions bounded away from zero or require
a different analysis.  Third, a shifted Hill function removes the
singularity as well, so this advantage is over the unshifted and leaky
forms, not over the Hill family as a whole.

\section{Discussion and conclusion}
\label{sec:conclusion}

For the two-gene loop with one measured species, observability and
Luenberger observer design can be carried out explicitly, and each
conclusion rests on a specific property of the regulation function.
Global injectivity of the observability map needs only strict
monotonicity; the rank condition at a point needs a non-zero slope
there; a uniform bound on the sensitivity of state reconstruction
needs a positive lower slope bound on the interval of interest; and
global exponential convergence with noise attenuation needs a finite
upper slope bound together with feasibility of the four-vertex matrix
inequality.  All of these hold for the logistic, Hill, leaky-Hill and
shifted-Hill families alike when the Hill exponent exceeds one, and the
four-vertex inequality was feasible for every model examined whose
slopes are bounded (for unshifted and leaky Hill functions with exponent
below one only set-based certificates exist).  A
positive lower slope bound can lift the certified rate above the
degradation rate of the unmeasured species, up to a closed-form
ceiling, when the corresponding inequality is feasible.  The logistic
function gives these constants in closed form and has a slope bounded
away from zero at the basal end; in exchange, its slope decays
exponentially towards saturation.  At matched threshold and slope, and
for Hill exponents above one, the logistic form is better conditioned
than leaky or shifted Hill forms when the threshold lies in the middle
of the expression range, and worse when it lies low in that range,
while certified convergence rates differ little between families.  For
Hill exponents below one the picture changes: the unshifted Hill slope
is unbounded at zero, so its observability map stays injective and its
inverse has a small Lipschitz constant at the basal end, but the
global bounded-slope certificate is unavailable there, while it
remains available for the logistic and shifted Hill functions.

The analysis has several limitations.  The matching at one point is a
convention, not a fit: away from the
threshold the curves differ, and so do the equilibria, bifurcation
structure and time courses of the corresponding models.  Whether a
logistic, leaky-Hill or shifted-Hill form describes a given regulatory
interaction better is an empirical question that requires
dose--response data and identifiability analysis, which we do not
address.  The observer~\eqref{eq:observer} assumes that the model and
its parameters are known; robustness to parameter uncertainty is not
analysed.  The extension of Hill-type functions to negative arguments
is needed only because a Luenberger observer does not preserve
positivity; positive observers~\citep{proskurnikov2019positive} avoid
it by construction, and the logistic function avoids this technical
step.

Natural extensions are the design of gains directly from the vertex
inequalities, robustness to parameter uncertainty, observers for
longer cascades (for which \ref{app:chain} gives the
observability structure), and a comparison based on models fitted to
the same dose--response data.

\section*{Declaration of competing interest}

The author declares that there are no known competing financial
interests or personal relationships that could have appeared to
influence the work reported in this paper.

\section*{Funding}

This research did not receive any specific grant from funding agencies
in the public, commercial, or not-for-profit sectors.

\section*{CRediT authorship contribution statement}

\textbf{Ismail Belgacem:} Conceptualization, Methodology, Formal
analysis, Software, Validation, Visualization, Writing -- original
draft, Writing -- review \& editing.

\section*{Data and code availability}

No experimental data were used.  The numerical results were obtained
with R~4.3.3 and the package \texttt{deSolve}~1.40.  The scripts used
for the simulations, the operating sets, the gains and the
verification of all matrix inequalities, together with the data file
listing every certificate of \ref{app:certificates}, are available
from the author upon request.

\appendix

\section{Observability of a regulatory cascade}
\label{app:chain}

Consider the cascade of $N$ genes
$\dot x_i=\kappa_ir_i(x_{i+1})-\gamma_ix_i$ ($i=1,\dots,N-1$),
$\dot x_N=\kappa_Nu-\gamma_Nx_N$ with a known constant input $u$,
$C^{N-1}$ regulation functions and $y=x_1$.

\begin{theorem}[Observability matrix of the cascade]
\label{thm:chain}
The Lie-derivative observability matrix $\mathcal{O}_N$ of the cascade
is lower triangular with diagonal entries $d_0=1$ and, for
$k=1,\dots,N-1$, $d_k=\prod_{i=1}^{k}\kappa_ir_i'(x_{i+1})$, so that
\begin{equation}
  \det\mathcal{O}_N(\mathbf{x})
  =\prod_{i=1}^{N-1}\bigl(\kappa_ir_i'(x_{i+1})\bigr)^{N-i}.
  \label{eq:ngene_det}
\end{equation}
Let $B=I_1\times\dots\times I_N$ be a box of intervals such that
$|r_i'|\ge m_i>0$ on $I_{i+1}$ for $i=1,\dots,N-1$.  Then the map
$\mathbf{x}\mapsto(y,\dot y,\dots,y^{(N-1)})$ is injective on $B$ and
$|d_k|\ge\prod_{i=1}^{k}\kappa_im_i$ on $B$.
\end{theorem}

\begin{proof}
$L_f^kh$ depends on $x_1,\dots,x_{k+1}$ only, and
$\partial L_f^kh/\partial x_{k+1}=d_k$ by induction: the only term of
$L_f^{k+1}h=\sum_j(\partial L_f^kh/\partial x_j)\dot x_j$ that depends
on $x_{k+2}$ is $d_k\kappa_{k+1}r_{k+1}(x_{k+2})$.  Each factor
$\kappa_ir_i'$ appears in $d_k$ for $k\ge i$, i.e.\ $N-i$ times.
Injectivity follows by solving successively for $x_1,x_2,\dots$, since
$L_f^kh$ is strictly monotone in $x_{k+1}$ for fixed
$x_1,\dots,x_k$ (it equals $d_{k-1}\kappa_kr_k(x_{k+1})$ plus a
function of $x_1,\dots,x_k$, and $d_{k-1}\ne0$ has constant sign).
\end{proof}

The product structure shows how small slopes compound along the
cascade; it applies to every family, with the family entering only
through the bounds $m_i$.  For quantitative statements the outputs
are scaled as in Section~\ref{sec:observability},
$(y,\tau\dot y,\dots,\tau^{N-1}y^{(N-1)})$, which multiplies $d_k$ by
$\tau^k$.  For $N=2$ the result reduces to the rank statement of
Theorem~\ref{thm:obs}(ii).

\section{Certificates for the reported rates}
\label{app:certificates}
\setcounter{table}{0}

Table~\ref{tab:certs} lists, for every rate quoted in
Section~\ref{sec:design} and Tables~\ref{tab:compare}
and~\ref{tab:subunit}, the certified bracket for $\alpha^\star/2$ and the
data that establish it.  The columns ``lower'' and ``upper'' are bounds
on the rate $\alpha^\star/2$, so the primal residual for a lower entry
$r_{\rm L}$ is evaluated at $\alpha=2r_{\rm L}$ and the dual matrix $G$
for an upper entry $r_{\rm U}$ at $\alpha=2r_{\rm U}$.  The lower bound
is the reported rate $\alpha_{\mathrm{cert}}/2$; it is certified by the
matrix $P$ shown, normalised to $\det P=1$ before rounding and rounded to
at least four significant digits (trailing zeros omitted), more where
needed: the column ``residual'' is
$\max_v\lambda_{\max}(\tilde A(v)^TP+P\tilde A(v)+\alpha P)$ evaluated
with the printed entries, not with the unrounded optimiser output.  The
upper bound is the lower bound plus $10^{-4}$; it follows either from
Proposition~\ref{prop:limit} (when $\bar\eta$ does not exceed it) or from a
dual certificate (Lemma~\ref{lem:dual}).  For the latter the column
$\lambda_{\min}(G)$ is computed from the factors $L_v$ rounded to the
fewest significant digits (between four and eight) that keep it
positive; the smallest margin, $4\times10^{-9}$, exceeds the
double-precision error of these $2\times2$ computations by about seven
orders of magnitude.  The factors $L_v$ of every dual certificate, and
all matrices of Tables~\ref{tab:certs}, \ref{tab:frho} and~\ref{tab:thr}, are listed in
the data file \texttt{certificates\_data.txt}, available with the
scripts.

Table~\ref{tab:frho} gives the finite-$\rho$ certificates
($\rho=0.5$~nM) quoted in Section~\ref{sec:subunit}; for these models
$r_1$ and $r_2$ have the same slope range on the enlarged interval, and
the admissible initial error is $\rho/\sqrt{\kappa(P)}$
(condition~\eqref{eq:rho_cond} with $W=0$), rounded down.
Table~\ref{tab:thr} gives the certificates on both sides of the
feasibility transitions of Section~\ref{sec:subunit}; feasibility at
$\varepsilon=0.368$ and $0.392$ is checked at half the computed optimal
rate.

\begin{table*}[t]
\centering
\caption{Certificates for the reported rates (min$^{-1}$).  Labels A, B,
C refer to the scenarios of Sections~\ref{sec:compare}
and~\ref{sec:subunit}; ``set'' is the operating-set (Table~\ref{tab:compare})
or interval (Table~\ref{tab:subunit}) certificate with $\rho=0$.}
\label{tab:certs}
\footnotesize
\begin{tabular}{@{}lccllllll@{}}
\toprule
Case & lower & upper & by & $p_{11}$ & $p_{12}$ & $p_{22}$ & residual & $\lambda_{\min}(G)$ \\
\midrule
Baseline, global & $0.0999$ & $0.1000$ & Prop.~3 & $6.994$ & $0.2648$ & $0.153$ & $-1.1\times10^{-6}$ & -- \\
Baseline, set, $\rho=0$ & $0.1111$ & $0.1112$ & Prop.~3 & $6.111$ & $0.3235$ & $0.1808$ & $-1.3\times10^{-5}$ & -- \\
Baseline, set, $\rho=1$ & $0.1041$ & $0.1042$ & Prop.~3 & $5.971$ & $0.3141$ & $0.184$ & $-3.6\times10^{-6}$ & -- \\
Baseline, set, $\rho=5$ & $0.1000$ & $0.1001$ & Prop.~3 & $6.9919$ & $0.2647$ & $0.15304$ & $-6.0\times10^{-6}$ & -- \\
\midrule
A: Logistic, set & $0.1119$ & $0.1120$ & Prop.~3 & $5.894$ & $0.3382$ & $0.1891$ & $-2.3\times10^{-5}$ & -- \\
A: Logistic, global & $0.0999$ & $0.1000$ & Prop.~3 & $6.994$ & $0.2648$ & $0.153$ & $-1.1\times10^{-6}$ & -- \\
A: Hill, set & $0.0999$ & $0.1000$ & Prop.~3 & $6.858$ & $0.2714$ & $0.1566$ & $-1.3\times10^{-5}$ & -- \\
A: Hill, global & $0.0999$ & $0.1000$ & Prop.~3 & $6.858$ & $0.2714$ & $0.1566$ & $-1.3\times10^{-5}$ & -- \\
A: Leaky Hill, set & $0.1001$ & $0.1002$ & Prop.~3 & $7.346$ & $0.2761$ & $0.1465$ & $-1.0\times10^{-5}$ & -- \\
A: Leaky Hill, global & $0.0999$ & $0.1000$ & Prop.~3 & $7.35$ & $0.2763$ & $0.1465$ & $-1.6\times10^{-5}$ & -- \\
A: Shifted Hill, set & $0.1001$ & $0.1002$ & dual & $6.895$ & $0.2691$ & $0.1555$ & $-1.1\times10^{-6}$ & $1.4\times10^{-6}$ \\
A: Shifted Hill, global & $0.0999$ & $0.1000$ & Prop.~3 & $6.899$ & $0.2693$ & $0.1555$ & $-1.4\times10^{-5}$ & -- \\
\midrule
B: Logistic, set & $0.0994$ & $0.0995$ & dual & $4.771$ & $0.398$ & $0.2428$ & $-2.1\times10^{-5}$ & $1.1\times10^{-6}$ \\
B: Logistic, global & $0.0994$ & $0.0995$ & dual & $4.771$ & $0.398$ & $0.2428$ & $-2.1\times10^{-5}$ & $1.1\times10^{-6}$ \\
B: Hill, set & $0.0994$ & $0.0995$ & dual & $5.582$ & $0.5498$ & $0.2333$ & $-3.0\times10^{-5}$ & $2.8\times10^{-6}$ \\
B: Hill, global & $0.0994$ & $0.0995$ & dual & $5.582$ & $0.5498$ & $0.2333$ & $-3.0\times10^{-5}$ & $2.8\times10^{-6}$ \\
B: Leaky Hill, set & $0.0999$ & $0.1000$ & dual & $7.11$ & $0.969$ & $0.2727$ & $-5.3\times10^{-6}$ & $8.5\times10^{-8}$ \\
B: Leaky Hill, global & $0.0978$ & $0.0979$ & dual & $7.1183$ & $0.9702$ & $0.27272$ & $-4.5\times10^{-5}$ & $3.8\times10^{-7}$ \\
B: Shifted Hill, set & $0.0999$ & $0.1000$ & dual & $5.333$ & $0.4821$ & $0.2311$ & $-3.7\times10^{-5}$ & $1.5\times10^{-6}$ \\
B: Shifted Hill, global & $0.0995$ & $0.0996$ & dual & $5.336$ & $0.48236$ & $0.23101$ & $-4.6\times10^{-6}$ & $7.4\times10^{-6}$ \\
\midrule
C: Logistic, global & $0.0999$ & $0.1000$ & Prop.~3 & $71.98$ & $0.2679$ & $0.01489$ & $-2.6\times10^{-6}$ & -- \\
C: Logistic, set $[0,10]$ & $0.7794$ & $0.7795$ & Prop.~3 & $166.5$ & $4.092$ & $0.1066$ & $-1.5\times10^{-6}$ & -- \\
C: Hill, set $[1,10]$ & $0.2741$ & $0.2742$ & dual & $30.67$ & $0.8511$ & $0.05622$ & $-2.0\times10^{-6}$ & $2.4\times10^{-7}$ \\
C: Leaky Hill, set $[1,10]$ & $0.2676$ & $0.2677$ & dual & $32.57$ & $0.8111$ & $0.0509$ & $-7.5\times10^{-6}$ & $3.3\times10^{-8}$ \\
C: Shifted Hill, global & $0.0929$ & $0.0930$ & dual & $31.87$ & $0.8843$ & $0.05592$ & $-1.2\times10^{-5}$ & $4.2\times10^{-9}$ \\
C: Shifted Hill, set $[0,10]$ & $0.2698$ & $0.2699$ & dual & $29.01$ & $0.8049$ & $0.0568$ & $-3.1\times10^{-6}$ & $4.6\times10^{-8}$ \\
\bottomrule
\end{tabular}
\end{table*}

\begin{table*}[t]
\centering
\caption{Finite-$\rho$ certificates for scenario C with $\rho=0.5$~nM
(Section~\ref{sec:subunit}).  Slope range of $|r_1'|$ (equal to that of
$r_2'$) on the enlarged interval (nM$^{-1}$), certified rate
(min$^{-1}$), matrix $P$ ($\det P=1$, as printed), residual at
$\alpha=2\alpha_{\mathrm{cert}}/2$ evaluated with the printed entries,
$\lambda_{\min}(P)$, $\sqrt{\kappa(P)}$ and the admissible initial error
$\rho/\sqrt{\kappa(P)}$ (nM, $W=0$, rounded down).}
\label{tab:frho}
\footnotesize
\begin{tabular}{@{}lllllllllll@{}}
\toprule
Family & $I^\rho$ & slopes & $\alpha_{\mathrm{cert}}/2$ & $p_{11}$ & $p_{12}$ & $p_{22}$ & residual & $\lambda_{\min}(P)$ & $\sqrt{\kappa(P)}$ & $\|\mathbf e(0)\|_{\max}$ \\
\midrule
Logistic & $[-0.5,10.5]$ & $[0.0232,0.0250]$ & $0.7584$ & $137.6$ & $3.594$ & $0.1011$ & $-1.6\times10^{-6}$ & $0.00722$ & $138.1$ & $0.0036$ \\
Hill & $[0.5,10.5]$ & $[0.0115,0.1825]$ & $0.0911$ & $21.226$ & $0.5883$ & $0.063418$ & $-4.5\times10^{-6}$ & $0.0471$ & $21.2$ & $0.0235$ \\
Leaky Hill & $[0.5,10.5]$ & $[0.0109,0.1734]$ & $0.0738$ & $22.77$ & $0.5665$ & $0.05801$ & $-7.6\times10^{-6}$ & $0.0439$ & $22.8$ & $0.0219$ \\
Shifted Hill & $[-0.5,10.5]$ & $[0.0126,0.1138]$ & $0.2621$ & $29.13$ & $0.8083$ & $0.05675$ & $-9.5\times10^{-6}$ & $0.0343$ & $29.2$ & $0.0171$ \\
\bottomrule
\end{tabular}
\end{table*}

\begin{table}[ht]
\centering
\caption{Certificates on both sides of the feasibility transitions of
Section~\ref{sec:subunit} (scenario C).}
\label{tab:thr}
\footnotesize
\begin{tabular}{@{}lll>{\raggedright\arraybackslash}p{0.44\linewidth}@{}}
\toprule
Model & Interval & Inequality & Certificate \\
\midrule
Hill & $[0.367,10]$ & infeasible & dual certificate at $\alpha=0$, $\lambda_{\min}(G)=5.2\times10^{-6}$ \\
Hill & $[0.368,10]$ & feasible & $P=\bigl(\begin{smallmatrix}18.478&0.51177\\0.51177&0.068293\end{smallmatrix}\bigr)$, residual $-2.9\times10^{-5}$ at $\alpha=3.2\times10^{-4}$ \\
Leaky Hill & $[0.391,10]$ & infeasible & dual certificate at $\alpha=0$, $\lambda_{\min}(G)=1.9\times10^{-6}$ \\
Leaky Hill & $[0.392,10]$ & feasible & $P=\bigl(\begin{smallmatrix}20.406&0.5074\\0.5074&0.061623\end{smallmatrix}\bigr)$, residual $-5.3\times10^{-5}$ at $\alpha=6.0\times10^{-4}$ \\
\bottomrule
\end{tabular}
\end{table}

\bibliographystyle{elsarticle-num}
\bibliography{refs}

\end{document}